\documentclass[a4paper]{amsart}
\usepackage{amsfonts}
\usepackage{amssymb}
\usepackage{amsmath}  %Overlap with newcommand \mod
\usepackage{verbatim}
\usepackage{stmaryrd}

\input xy
\xyoption{all}

\newtheorem{theorem}{Theorem}[section]
\newtheorem{lemma}[theorem]{Lemma}

\newtheorem{proposition}[theorem]{Proposition}
\newtheorem{corollary}[theorem]{Corollary}
\newtheorem{definition}[theorem]{Definition}
\newtheorem{remark}[theorem]{Remark}

\newcommand{\Ker}{\mbox{\rm Ker\,}}

\newcommand{\Ima}{\mbox{\rm Im\,}}
\newcommand{\End}[1]{\mbox{\rm End}_{#1}}
\newcommand{\Aut}[1]{\mbox{\rm Aut}_{#1}}
\newcommand{\Hom}[1]{\mbox{{\rm Hom}}_{#1}}
\newcommand{\Ext}[2]{\mbox{\rm Ext}^{#1}_{#2}}

\newcommand{\pr}[1]{\mbox{{\rm pr}}_{#1}}
\newcommand{\add}{\mbox{{\rm add \!}}}
\newcommand{\MOD}{\mbox{{\rm mod \!}}}
\newcommand{\Mod}{\mbox{{\rm Mod \!}}}

\newcommand{\perf}{\mbox{{\rm per \!}}}
\newcommand{\Ho}{\mbox{{\rm H}}}
\newcommand{\ind}[1]{\mathrm{ind}_{#1}}

\newcommand{\Gr}[1]{\mbox{{\rm Gr}}_{#1}}

\newcommand{\K}{\mbox{{\rm K}}}

\newcommand{\Lim}{\mbox{{\rm lim \!}}}

\newcommand{\demo}[1]{\textsc{Proof} #1 \hfill $\Box$ \bigskip}

\newcommand{\cA}{\mathcal{A}}

\newcommand{\cC}{\mathcal{C}}
\newcommand{\cD}{\mathcal{D}}

\newcommand{\cF}{\mathcal{F}}
\newcommand{\cH}{\mathcal{H}}
\newcommand{\cM}{\mathcal{M}}
\newcommand{\cN}{\mathcal{N}}

\newcommand{\cT}{\mathcal{T}}
\newcommand{\cU}{\mathcal{U}}
\newcommand{\bN}{\mathbb{N}}
\newcommand{\bP}{\mathbb{P}}
\newcommand{\bQ}{\mathbb{Q}}
\newcommand{\bZ}{\mathbb{Z}}
\newcommand{\bg}{\mathbf{g}}
\newcommand{\bh}{\mathbf{h}}
\newcommand{\bx}{\mathbf{x}}
\newcommand{\by}{\mathbf{y}}
\newcommand{\fm}{\mathfrak{m}}

\begin{document}

\title[Cluster algebras via cluster categories]{Cluster algebras via cluster categories with infinite-dimensional morphism spaces}
\author{Pierre-Guy Plamondon}
\address{Universit\'e Paris Diderot -- Paris 7\\
   Institut de Math\'ematiques de Jussieu, UMR 7586 du CNRS \\
   Case 7012\\
   B\^atiment Chevaleret\\
   75205 Paris Cedex 13\\
   France }
\email{plamondon@math.jussieu.fr}
\thanks{The author was financially supported by an NSERC scholarship.}
%\date{\today}

\begin{abstract}
We apply our previous work on cluster characters for $\Hom{}$-infinite cluster categories to the theory of cluster algebras.  We give a new proof of Conjectures 5.4, 6.13, 7.2, 7.10 and 7.12 of Fomin and Zelevinsky's \emph{Cluster algebras~IV} for skew-symmetric cluster algebras.  We also construct an explicit bijection sending certain objects of the cluster category to the decorated representations of Derksen, Weyman and Zelevinsky, and show that it is compatible with mutations in both settings.  Using this map, we give a categorical interpretation of the $E$-invariant and show that an arbitrary decorated representation with vanishing $E$-invariant is characterized by its $\bg$-vector.  Finally, we obtain a substitution formula for cluster characters of not necessarily rigid objects.
\end{abstract}

\maketitle

\tableofcontents

%========================================================================================================
\section{Introduction}
Since their introduction by S.~Fomin and A.~Zelevinsky in \cite{FZ02}, cluster algebras have been found to enjoy connections with several branches of mathematics, see for instance the survey papers \cite{Zelevinsky07}, \cite{GLS07} and \cite{Keller09}.

Cluster algebras are commutative algebras generated by \emph{cluster variables} grouped into sets of fixed finite cardinality called \emph{clusters}.  Of particular importance are cluster algebras \emph{with coefficients}, as most known examples of cluster algebras are of this kind.  In this paper, we will work with cluster algebras of geometric type with coefficients.

In \cite{FZ07}, the authors developped a combinatorial framework allowing  the study of coefficients in cluster algebras.  Important tools that the authors introduced are the $F$-polynomials and $\bg$-vectors.  In particular, they proved that the behaviour of the coefficients in any cluster algebra is governed by that of the coefficients in a cluster algebra \emph{with principal coefficients}, using the $F$-polynomials (see \cite[Theorem 3.7]{FZ07}).

The authors phrased a number of conjectures, mostly regarding $F$-polynomials and $\bg$-vectors.  We list some of them here:
\begin{description}
	\item[(5.4)] every $F$-polynomial has constant term $1$;
	\item[(6.13)] the $\bg$-vectors of the cluster variables of any given seed are \emph{sign-coherent} in a sense to be defined;
	\item[(7.2)] cluster monomials are linearly independent;
	\item[(7.10)] different cluster monomials have different $\bg$-vectors, and the $\bg$-vectors of the cluster variables of any cluster form a basis of $\bZ^r$;
	\item[(7.12)] the mutation rule for $\bg$-vectors can be expressed using a certain piecewise-linear transformation. 
\end{description}

Work on these conjectures includes
\begin{itemize}
	\item a proof of (7.2) by P.~Sherman and A.~Zelevinsky \cite{SZ04} for Dynkin and affine type of rank $2$;
	\item a proof of (7.2) by P.~Caldero and B.~Keller \cite{CK08} for Dynkin type;
	\item a proof of (7.2) by G.~Dupont \cite{Dupont08} for affine type $A$;
	\item a proof of (7.2) by M.~Ding, J.~Xiao and F.~Xu \cite{DXX09} for affine types;
	\item a proof of (7.2) by G.~Cerulli Irelli \cite{Cerulli09} in type $A_2^{(1)}$ by explicit computations;
	\item a proof of (5.4) by R.~Schiffler \cite{Schiffler10} for cluster algebras arising from unpunctured surfaces;
	\item a proof of (7.2) by L.~Demonet \cite{Demonet09} for certain skew-symmetrizable cluster algebras;
	\item a proof of all five conjectures by C.~Fu and B.~Keller \cite{FK09} for cluster algebras categorified by $\Hom{}$-finite $2$-Calabi--Yau Frobenius or triangulated categories, using work of R.~Dehy and B.~Keller \cite{DK08};
	\item a proof of (7.2) by C.~Geiss, B.~Leclerc and J.~Schr\"oer \cite{GLS08} for acyclic cluster algebras;
	\item a proof of (5.4), (6.13), (7.10) and (7.12) in full generality by H.~Derksen, J.~Weyman and A.~Zelevinsky \cite{DWZ09} using decorated representations of quivers with potentials;
	\item a proof of (5.4), (6.13), (7.10) and (7.12) in full generality by K.~Nagao \cite{Nag09} using Donaldson--Thomas theory (see for instance \cite{JS08}, \cite{KS08} and \cite{Bridgeland10}).
\end{itemize}

In this paper, we use (generalized) cluster categories to give a new proof of (5.4), (6.13), (7.10) and (7.12) in full generality, and to prove (7.2) for any skew-symmetric cluster algebra of geometric type whose associated matrix is of full rank.

More precisely, we use the cluster category introduced by A.~Buan, R.~Marsh, M.~Reineke, I.~Reiten and G.~Todorov in \cite{BMRRT06} (and independently by P.~Caldero, F.~Chapoton and R.~Schiffler in \cite{CCS06} in the $A_n$ case) and generalized to any quiver with potential by C.~Amiot in \cite{Amiot08}.  Note that this category can be $\Hom{}$-infinite.  We obtain applications to cluster algebras via the cluster character of Y.~Palu \cite{Palu08}, which generalized the map introduced by P.~Caldero and F.~Chapoton in \cite{CC06}.  It was extended in \cite{Plamondon09} to possibly $\Hom{}$-infinite cluster categories.  In particular, we have to restrict the cluster character to a suitable subcategory $\cD$ of the cluster category.

Using this cluster character, we give categorical interpretations of $F$-polynomials and $\bg$-vectors which allow us to prove the conjectures mentioned above.  We prove (7.2), (6.13), (7.10) and (7.12) in section \ref{sect::index g-vectors} and (5.4) in section \ref{sect::clustercharF}. Some of our results concerning rigid objects in section \ref{sect::rigid} and indices in section \ref{sect::index g-vectors} are used in \cite{IIKKN1} and \cite{IIKKN2}.  The methods we use are mainly generalizations of those used for the $\Hom{}$-finite case in \cite{DK08} and \cite{FK09}. 

The key tool that we use in our proofs is the multiplication formula proved in \cite[Proposition 3.16]{Plamondon09}, while the proofs of H.~Derksen, J.~Weyman and A.~Zelevinsky rely on a substitution formula \cite[Lemma 5.2]{DWZ09}.

We also show in section \ref{sect::dwz} that the setup used in \cite{DWZ09} is closely related to the cluster-categorical approach.  We prove in section \ref{sect::bijection} that (isomorphism classes of) decorated representations over a quiver with potential are in bijection with  (isomorphism classes of) objects in the subcategory $\cD$ of the cluster category.  In sections \ref{sect::fg} and \ref{sect::E}, we give  an interpretation of the $F$-polynomial, $\bg$-vector, $h$-vector and $E$-invariant of a decorated representation in cluster-categorical terms.  In particular, we prove a stronger version of \cite[Lemma 9.2]{DWZ09} in Corollary \ref{coro::dwz}.  Using the substitution formula for $F$-polynomials \cite[Lemma 5.2]{DWZ09}, we also obtain a substitution formula for cluster characters of not necessarily rigid object (Corollary \ref{coro::subst}).

\subsection*{Acknowledgements}
This work is part of my PhD thesis, supervised by Professor Bernhard Keller.  I would like to thank him here for his generosity both in time and mathematical knowledge.

%========================================================================================================
\section{Recollections}\label{sect::recollections}

%--------------------------------------------------------------------------------------
\subsection{Background on cluster algebras}
We give a brief summary of the definitions and results we will need concerning cluster algebras.  Our main source for the material in this section is \cite{FZ07}.

%......................................................
\subsubsection{Cluster algebras with coefficients}

We will first recall the definition of (skew-symmetric) cluster algebras (of geometric type).

The \emph{tropical semifield} $\textrm{Trop}(u_1, u_2, \ldots, u_n)$ is the abelian group (written multiplicatively) freely generated by the $u_i$'s, with an \emph{addition} $\oplus$ defined by
\begin{displaymath}
	\prod_j u_j^{a_j} \oplus \prod_j u_j^{b_j} = \prod_j u_j^{\textrm{min}(a_j, b_j)}.
\end{displaymath}

A \emph{quiver} is an oriented graph.  Thus, it is given by a quadruple $Q = (Q_0, Q_1, s, t)$, where $Q_0$ is the set of \emph{vertices}, $Q_1$ is the set of \emph{arrows}, and $s$ (respectively $t$) is a map from $Q_1$ to $Q_0$ which sends each arrow to its \emph{source} (respectively its \emph{target}).  A quiver is \emph{finite} if it has finitely many vertices and arrows.

An \emph{ice quiver} (see \cite{FK09}) is a pair $(Q,F)$, where $Q$ is a quiver and $F$ is a subset of $Q_0$.  The elements of $F$ are the \emph{frozen vertices} of $Q$.  It is \emph{finite} if $Q$ is finite.

Let $r$ and $n$ be integers such that $1\leq r \leq n$.  Denote by $\bP$ the tropical semifield $\textrm{Trop}(x_{r+1}, \ldots, x_n)$.  Let $\cF$ be the field of fractions of the ring of polynomials in $r$ indeterminates with coefficients in $\bQ\bP$.

Let $(Q,F)$ be a finite ice quiver, where $Q$ has no oriented cycles of length $\leq 2$, and $F$ and $Q_0$ have $r$ and $n$ elements respectively.  We will denote the elements of $Q_0\setminus F$ by the numbers $1, 2, \ldots, r$  and those of $F$ by $(r+1), (r+2), \ldots, n$.  Let $i$ be in $Q_0 \setminus F$.  One defines the \emph{mutation of $(Q,F)$ at $i$} as the ice quiver $\mu_i(Q,F) = (Q',F')$ constructed from $(Q,F)$ as follows:
\begin{itemize}
	\item the sets $Q'_0$ and $F'$ are equal to $Q_0$ and $F$, respectively;
	
	\item all arrows not incident with $i$ in $Q$ are kept;
	
	\item for each subquiver of $Q$ of the form $j\rightarrow i \rightarrow \ell$, an arrow from $j$ to $\ell$ is added;
	
	\item all arrows incident to $i$ are reversed;
	
	\item arrows from a maximal set of pairwise disjoint oriented cycles of length two in the resulting quiver are removed.
	
\end{itemize}

A \emph{seed} is a pair $\big( (Q,F), \bx  \big)$, where $(Q,F)$ is an ice quiver as above, and $\bx = \{x_1, \ldots, x_r \}$ is a free generating set of the field $\cF$.  Given a vertex $i$ of $Q_0 \setminus F$, the \emph{mutation} of the seed $\big( (Q,F), \bx \big)$ at the vertex $i$ is the pair $\mu_i\big( (Q,F), \bx \big) = \big( (Q', F'), \bx' \big)$, where
\begin{itemize}
	\item $(Q', F')$ is the mutated ice quiver $\mu_i(Q, F)$;
	\item $\bx' = \bx \setminus \{x_i\} \cup \{x'_i\}$, where $x'_i$ is obtained from the \emph{exchange relation}
	   \begin{displaymath}
	      x_ix'_i = \prod_{\substack{ \alpha \in Q_1 \\ s(\alpha) = i}} x_{t(\alpha)} +  \prod_{\substack{ \alpha \in Q_1 \\ t(\alpha) = i}} x_{s(\alpha)}.
     \end{displaymath}
\end{itemize}
The mutation of a seed is still a seed, and the mutation at a fixed vertex is an involution.

Fix an \emph{initial seed} $\big( (Q,F), \bx \big)$.
\begin{itemize}
	\item The sets $\bx'$ obtained by repeated mutation of the initial seed are the \emph{clusters}.
	\item The elements of the clusters are the \emph{cluster variables}.
	\item The $\bZ\bP$-subalgebra of $\cF$ generated by all cluster variables is the \emph{cluster algebra} $\cA = \cA\big( (Q,F), \bx \big)$.
\end{itemize} 

Suppose that $n=2r$.  A cluster algebra has \emph{principal coefficients at a seed $\big( (Q',F'), \bx' \big)$} if there is exactly one arrow from $(r+\ell)$ to $\ell$ (for $\ell = 1,2,\ldots, r$), and if these are the only arrows whose source or target lies in $F'$.

%......................................................
\subsubsection{Cluster monomials and $\bg$-vectors}\label{sect::g-vectors}
Given an ice quiver $(Q,F)$, we associate to it an $(n\times r)$-matrix $\tilde{B} = (b_{ij})$, where each entry $b_{ij}$ is the number of arrows from $i$ to $j$ minus the number or arrows from $j$ to $i$.

Let $\big( (Q,F), \bx \big)$ be a seed of a cluster algebra $\cA$.  A \emph{cluster monomial} in $\cA$ is a product of cluster variables lying in the same cluster.

To define \emph{$\bg$-vectors}, we shall need a bit of notation.

For any integer $j$ between $1$ and $r$, let $\hat{y}_j$ be defined as
\begin{displaymath}
	\hat{y}_j = \prod_{\ell \in Q_0}x_{\ell}^{b_{\ell j}}.
\end{displaymath}

Let $\cM$ be the set of non-zero elements of $\cA$ which can be written in the form
\begin{displaymath}
	z = R(\hat{y}_1, \ldots, \hat{y}_r)\prod_{j=1}^{n}x_j^{g_j},
\end{displaymath}
where $R(u_1, \ldots, u_r)$ is an element of $\bQ(u_1, \ldots, u_r)$.  Note that all cluster monomials belong to $\cM$.

By \cite[Proposition 7.8]{FZ07}, if the matrix $\tilde{B}$ is of full rank $r$, then any element of $\cM$ can be written in a unique way in the form above, with $R$ \emph{primitive} (that is, $R$ can be written as a ration of two polynomials, none of which is divisible by any of the $u_j$'s).  In that case, if $z$ is an element of $\cM$ written as above with $R$ primitive, the vector
\begin{displaymath}
	\bg(z) = (g_1, \ldots, g_r)
\end{displaymath}
is the \emph{$\bg$-vector} of $z$.

Let us now state Conjectures 7.2, 7.10 and 7.12 of \cite{FZ07}.

\begin{description}
	\item[7.2] Cluster monomials are linearly independent over $\bZ\bP$.
	\item[7.10] Different cluster monomials have different $\bg$-vectors; the $\bg$-vectors of the cluster variables of any cluster form a $\bZ$-basis of $\bZ^r$.
	\item[7.12] Let $\bg = (g_1, \ldots, g_r)$ and $\bg' = (g'_1, \ldots, g'_r)$ be the $\bg$-vectors of one cluster monomial with respect to two clusters $t$ and $t'$ related by one mutation at the vertex $i$.  Then we have
	\begin{displaymath}
	  g'_j = \left\{ \begin{array}{ll}
                      -g_i & \textrm{if } j=i\\
                      g_j + [b_{ji}]_{+}g_{i} - b_{ji}\min(g_{i}, 0)  & \textrm{if } j\neq i
                   \end{array} \right.
  \end{displaymath}
  where $B = (b_{j\ell})$ is the matrix associated with the seed $t$, and we set $[x]_{+} = \max(x,0)$ for any real number $x$.
\end{description}

%......................................................
\subsubsection{$F$-polynomials}\label{sect::F-polynomials}

Let $\cA$ be a cluster algebra with principal coefficients at a seed $\big( (Q,F), \bx \big)$.  Let $t$ be a seed of $\cA$ and $\ell$ be a vertex of $Q$ that is not in $F$.  Then the $\ell$-th cluster variable of $t$ can be written as a subtraction-free rational function in variables $x_1, \ldots, x_{2r}$.  Following \cite[Definition 3.3]{FZ07}, we define the \emph{$F$-polynomial} $F_{\ell, t}$ as the specialization of this rational function at $x_1 = \ldots = x_r = 1$.  It was proved in \cite[Proposition 3.6]{FZ07} that $F_{\ell, t}$ is indeed a polynomial.

We now state Conjecture 5.4 of \cite{FZ07} : \emph{Every $F$-polynomial has constant term $1$}. 

%......................................................
\subsubsection{$Y$-seeds and their mutations}\label{sect::tropical}
We now recall  the notion of $Y$-seeds from \cite{FZ07}.  As above, let $1\leq r \leq n$ be integers, and let $\bP$ be the tropical semifield in the variables $x_{r+1}, \ldots, x_n$. 

A \emph{$Y$-seed} is a pair $(Q, \by)$, where
\begin{itemize}
	\item $Q$ is a finite quiver without oriented cycles of length $\leq 2$; and
	\item $\by = (y_1, \ldots, y_r)$ is an element of $\bP^r$.
\end{itemize}

Let $(Q, \by)$ be a $Y$-seed, and let $i$ be a vertex of $Q$.  The \emph{mutation of $(Q, \by)$ at the vertex $i$} is the $Y$-seed $\mu_i(Q, \by) = (Q', \by')$, where
\begin{itemize}
	\item $Q$ is the mutated quiver $\mu_i(Q)$; and
	\item $\by' = (y'_1, \ldots, y'_r)$ is obtained from $\by$ using the mutation rule
	   \begin{displaymath}
	       y'_j = \left\{ \begin{array}{ll}
              y_i^{-1} & \textrm{if $i=j$}\\
              y_j y_i^{m}(y_i \oplus 1)^{-m} & \textrm{if there are $m$ arrows from $i$ to $j$}\\
              y_j (y_i \oplus 1)^{m} & \textrm{if there are $m$ arrows from $j$ to $i$.}
                     \end{array} \right.
     \end{displaymath}
\end{itemize}

If, to any seed $\big( (Q,F), \bx\big)$ of a cluster algebra, we associate a $Y$-seed $(Q, \by)$ defined by
\begin{displaymath}
	y_j = \prod_{i=r+1}^n x_i^{b_{ij}},
\end{displaymath}
then for any such seed and its associated $Y$-seed, and for any vertex $i$ of $Q$, we have that the $Y$-seed associated to $\mu_i\big((Q,F), \bx\big)$ is $\mu_i(Q, \by)$.  This was proved in \cite{FZ07} after Definition 2.12.

%--------------------------------------------------------------------------------------
\subsection{Quivers with potentials and their mutations}\label{sect::QP}
We recall the notion of quiver with potential from \cite{DWZ08}.

Let $Q$ be a finite quiver.  Denote by $\widehat{kQ}$ its \emph{completed path algebra}, that is, the $k$-algebra whose underlying $k$-vector space is 
\begin{displaymath}
	\prod_{w \textrm{ path}}kw
\end{displaymath}
and whose multiplication is deduced from the composition of paths by distributivity (by convention, we compose paths from right to left).  It is a topological algebra for the $\fm$-adic topology, where $\fm$ is the ideal of $\widehat{kQ}$ generated by the arrows of $Q$.

A \emph{potential} on $Q$ is an element $W$ of the space
\begin{displaymath}
	Pot(Q) = \widehat{kQ}/C,
\end{displaymath}
where $C$ is the closure of the commutator subspace $[\widehat{kQ},\widehat{kQ}]$ in $\widehat{kQ}$.  In other words, it is a (possibly infinite) linear combination of cyclically inequivalent oriented cycles of $Q$.  The pair $(Q,W)$ is a \emph{quiver with potential}.

Given any arrow $a$ of $Q$, the \emph{cyclic derivative with respect to $a$} is the linear map $\partial_a$ from $Pot(Q)$ to $\widehat{kQ}$ whose action on (equivalence classes of) oriented cycles is given by
\begin{displaymath}
	\partial_a(b_r\cdots b_2 b_1) = \sum_{b_i = a}b_{i-1}b_{i-2}\cdots b_1b_rb_{r-1}\cdots b_{i+1}. 
\end{displaymath}

The \emph{Jacobian algebra} $J(Q,W)$ of a quiver with potential $(Q,W)$ is the quotient of the algebra $\widehat{kQ}$ by the closure of the ideal generated by the cyclic derivatives $\partial_a W$, as $a$ ranges over all arrows of $Q$.

The above map is generalized as follows.  For any path $p$ of $Q$, define $\partial_p$ as the linear map from $Pot(Q)$ to $\widehat{kQ}$ whose action on any (equivalence class of) oriented cycle $c$ is given by
\begin{displaymath}
	\partial_p(c) = \sum_{c = upv}vu + \sum_{\substack{c = p_1wp_2 \\ p=p_2p_1}}w,
\end{displaymath}
where the sums are taken over all decompositions of $c$ into paths of smaller length, with $u$, $v$ and $w$ possibly trivial paths, and $p_1$ and $p_2$ non-trivial paths.

Let $(Q,W)$ be a quiver with potential.  In order to define the \emph{mutation of $(Q,W)$ at a vertex $\ell$}, we must recall the process of reduction  of a quiver with potential.

Let $R$ be the $k$-algebra given by $\bigoplus_{i\in Q_0}ke_i$, where $e_i$ is the idempotent associated with the vertex $i$.  Two quivers with potentials $(Q,W)$ and $(Q',W')$ are \emph{right-equivalent} if $Q_0 = Q'_0$ and there exists an $R$-algebra isomorphism $\varphi : \widehat{kQ}\longrightarrow \widehat{kQ'}$ sending the class of $W$ to the class of $W'$ in $Pot(Q')$.

A quiver with potential $(Q,W)$ is \emph{trivial} if $W$ is a (possibly infinite) linear combination of paths of length at least $2$, and $J(Q,W)$ is isomorphic to $R$.  It is \emph{reduced} if $W$ has no terms which are cycles of length at most $2$. 

\begin{theorem}[\cite{DWZ08}, Theorem 4.6 and Proposition 4.5]
Any quiver with potential $(Q,W)$ is right equivalent to a direct sum of a reduced one $(Q_{red}, W_{red})$ and a trivial one $(Q_{triv}, W_{triv})$, both unique up to right-equivalence.  Moreover, $J(Q,W)$ and $J(Q_{red}, W_{red})$ are isomorphic.
\end{theorem}

We can now define the mutation of quivers with potentials.  Let $(Q,W)$ be a quiver with potential, and let $\ell$ be a vertex of $Q$ not involved in any cycle of length $\leq 2$.  Assume that $W$ is written as a series of oriented cycles which do not begin or end in $\ell$ ($W$ is always cyclically equivalent to such a potential).  The \emph{mutation of $(Q,W)$ at vertex $\ell$} is the new quiver with potential $\mu_{\ell}(Q,W)$ obtained from $(Q,W)$ as follows.
\begin{enumerate}
  \item For any subquiver $\xymatrix{i \ar[r]^{a} & \ell \ar[r]^{b} & j}$ of $Q$, add an arrow $\xymatrix{i \ar[r]^{[ba]} & j}$.
  \item Delete any arrow $a$ incident with $\ell$ and replace it by an arrow $a^{\star}$ going in the opposite direction; the first two steps yield a new quiver $\widetilde{Q}$.
  \item Let $\widetilde{W}$ be the potential on $\widetilde{Q}$ defined by $\widetilde{W} = [W] + \sum a^{\star}b^{\star}[ba]$, where the sum is taken over all subquivers of $Q$ of the form $\xymatrix{i \ar[r]^{a} & \ell \ar[r]^{b} & j}$, and where $[W]$ is obtained from $W$ by replacing each occurence of $ba$ in its terms by $[ba]$.  The first three steps yield a new quiver with potential $\widetilde{\mu}_{\ell}(Q,W) = (\widetilde{Q}, \widetilde{W})$.   
\end{enumerate} 

The mutation $\mu_{\ell}(Q,W)$ is then defined as the reduced part of $\widetilde{\mu}_{\ell}(Q, W)$. 

Note that $\mu_{\ell}(Q, W)$ might contain oriented cycles of length $2$, even if $(Q,W)$ did not.  This prevents us from performing iterated mutations following an arbitrary sequence of vertices.  

A vertex $i$ of $(Q,W)$ which is not involved in any oriented cycle of length $\leq 2$ (and thus at which mutation can be performed) is an \emph{admissible vertex}. An \emph{admissible sequence of vertices} is a sequence $\underline{i} = (i_1, \ldots, i_s)$ of vertices of $Q$ such that $i_1$ is an admissible vertex of $(Q,W)$, and $i_m$ is an admissible vertex of $\mu_{m-1}\mu_{m-2}\cdots\mu_{1}(Q,W)$, for $1 < m \leq s$.  In that case, we denote by $\mu_{\underline{i}}(Q,W)$ the mutated quiver with potential $\mu_{s}\mu_{s-1}\cdots\mu_{1}(Q,W)$.

%--------------------------------------------------------------------------------------
\subsection{Decorated representations and their mutations}\label{sect::deco}
We now recall from \cite[Section 10]{DWZ08} the notion of decorated representation of a quiver with potential.

Let $(Q,W)$ be a quiver with potential, and let $J(Q,W)$ be its Jacobian algebra.  A \emph{decorated representation} of $(Q,W)$ is a pair $\cM = (M, V)$, where $M$ is a finite-dimensional nilpotent $J(Q,W)$-right module and $V$ is a finite-dimensional $R$-module (recall that $R$ is given by $\bigoplus_{i\in Q_0} ke_i$).

We now turn to the \emph{mutation} of decorated representations.

Given a decorated representation $\cM = (M,V)$ of $(Q,W)$, and given any admissible vertex $\ell$ of $(Q,W)$, we construct a decorated representation $\widetilde{\mu}_{\ell}(\cM) = (\overline{M}, \overline{V})$ of $\overline{\mu}_{\ell}(Q,W)$ as follows.

We view $M$ as a representation of the opposite quiver $Q^{op}$ (we must work over the opposite quiver, since we use \emph{right} modules).  In particular, to each vertex $i$ is associated a vector space $M_i$ (which is equal to $Me_i$), and to each arrow $a:i\rightarrow j$ is associated a linear map $M_a:M_j \rightarrow M_i$.  For any path $p=a_r\cdots a_2a_1$, we denote by $M_p$ the linear map $M_{a_1}  M_{a_2}  \cdots  M_{a_r}$, and for any (possibly infinite) linear combination $\sigma = \sum_{i\in I}\lambda_i p_i$ of paths, we denote by $M_\sigma$ the linear map $\sum_{i\in I}\lambda_i M_{p_i}$ (this sum is finite since $M$ is nilpotent).  If $\sigma$ is zero in $J(Q,W)$, then $M_\sigma$ is the zero map.

Define the vector spaces $M_{in}$ and $M_{out}$ by
\begin{displaymath}
	M_{in} = \bigoplus_{\substack{a\in Q_1 \\ s(a) = \ell}}M_{t(a)} \quad \textrm{and} \quad M_{out} = \bigoplus_{\substack{b\in Q_1 \\ t(b) = \ell}}M_{s(b)}.
\end{displaymath}

Define the linear map $\alpha:M_{in}\longrightarrow M_\ell$ as the map given in matrix form by the line vector $\big( M_a:M_{t(a)}\rightarrow M_\ell \ \big| \ a\in Q_1, \ s(a) = \ell \big)$.  Similarly, define $\beta:M_\ell\longrightarrow M_{out}$ as the map given in matrix form by the column vector $\big( M_b:M_{\ell}\rightarrow M_{s(b)} \ \big| \ b\in Q_1, \ t(b) = \ell \big)$.

Define a third map $\gamma:M_{out}\longrightarrow M_{in}$ as the map given in matrix form by 
\begin{displaymath}
	\big( M_{\partial_{ab}W}:M_{s(b)}\rightarrow M_{t(a)} \ \big| \ a,b\in Q_1, \ s(a) = t(b) = \ell   \big).
\end{displaymath}

Now construct $\widetilde{\mu}_\ell(\cM) = (\overline{M}, \overline{V})$ as follows.
\begin{itemize}
	\item For any vertex $i\neq \ell$, set $\overline{M}_i = M_i$ and $\overline{V}_i = V_i$.
	
	\item Define $\overline{M}_\ell$ and $\overline{V}_\ell$ by 
	       \begin{displaymath}
	         \overline{M}_\ell = \frac{\Ker \gamma}{\Ima \beta} \oplus \Ima \gamma \oplus \frac{\Ker \alpha}{\Ima \gamma} \oplus V_\ell \quad \textrm{and} \quad \overline{V}_\ell = \frac{\Ker\beta}{\Ker\beta \cap \Ima \alpha}.
         \end{displaymath}
  
  \item For any arrow $a$ of $Q$ not incident with $\ell$, set $\overline{M}_a = M_a$.
  
  \item For any subquiver of the form $\xymatrix{i\ar[r]^a & \ell\ar[r]^b & j}$, set $\overline{M}_{[ba]} = M_{ba}$.
  
  \item the actions of the remaining arrows are encoded in the maps
         \begin{displaymath}
	         \overline{\alpha} =
                            \left( \begin{array}{c}
                          -\pi\rho \\
                          -\gamma \\
                            0 \\
                            0
                                   \end{array} \right) \quad \textrm{and} \quad \overline{\beta} = \left( \begin{array}{cccc}
                                                                                                                0 & \iota & \iota\sigma & 0                                                                                                                      \end{array} \right),
         \end{displaymath}
         where
         \begin{itemize}
	           \item the map $\rho:M_{out}\rightarrow \Ker\gamma$ is such that its composition with the inclusion map of $\Ker \gamma$ gives the identity map of $\Ker\gamma$;
	           \item the map $\pi:\Ker \gamma \rightarrow \Ker \gamma/\Ima \beta$ is the natural projection map;
	           \item the map $\sigma:\Ker\alpha / \Ima\gamma \rightarrow \Ker \alpha$ is such that its composition with the projection map $\Ker \alpha \rightarrow \Ker\alpha/\Ima\gamma$ gives the identity map of $\Ker\alpha/\Ima \gamma$;
	           \item the map $\iota:\Ima\gamma\rightarrow M_{in}$ is the natural inclusion map.
         \end{itemize}
          
\end{itemize}

It is shown in \cite[Proposition 10.7]{DWZ08} that $\widetilde{\mu}_\ell(\cM)$ is a decorated representation of $\widetilde{\mu}_\ell(Q,W)$.

%--------------------------------------------------------------------------------------
\subsection{Some invariants of decorated representations}\label{sect::invariants}
In this section, we recall from \cite{DWZ08} and \cite{DWZ09} the definitions of $F$-polynomial, $\bg$-vector, $h$-vector and $E$-invariant of a decorated representation.

We fix a quiver with potential $(Q,W)$ and a decorated representation $\cM = (M,V)$ of $(Q,W)$.  We number the vertices of $Q$ from $1$ to $n$.

The \emph{$F$-polynomial} of $\cM$ is the polynomial of $\bZ[u_1, \ldots , u_n]$ defined by
\begin{displaymath}
	F_{\cM}(u_1, \ldots, u_n) = \sum_{e}\chi\big( \Gr{e}(M) \big) \prod_{i=1}^{n}u_i^{e_i}.
\end{displaymath}

The \emph{$\bg$-vector} of $\cM$ is the vector $\bg_{\cM} = (g_1, \ldots, g_n)$ of $\bZ^n$, where
\begin{displaymath}
	g_i = \dim\Ker\gamma_i - \dim M_i + \dim V_i,
\end{displaymath}
where $\gamma_i$ is the map $\gamma : M_{out}\longrightarrow M_{in}$ defined in section \ref{sect::deco}.

The \emph{$h$-vector} of $\cM$ is the vector $h_{\cM} = (h_1, \ldots, h_n)$ of $\bZ^n$, where
\begin{displaymath}
	h_i = - \dim \Ker \beta_i
\end{displaymath}
where $\beta_i$ is the map $\beta : M_i\longrightarrow M_{out}$ defined in section \ref{sect::deco}.

The \emph{$E$-invariant} of $\cM$ is the integer
\begin{displaymath}
	E(\cM) = \dim\Hom{J(Q,W)}(M, M) + \sum_{i=1}^n g_i \dim M_i,
\end{displaymath}
where $(g_1, \ldots, g_n)$ is the $\bg$-vector of $\cM$.

Let $\cN = (N, U)$ be another decorated representation of $(Q,W)$.  The $E$-invariant can also be defined using the two integer-valued invariants
\begin{displaymath}
	E^{inj}(\cM, \cN) = \dim\Hom{J(Q,W)}(M, N) + \sum_{i=1}^n (\dim M_i)g_i(\cN) \quad \textrm{and}
\end{displaymath}
\begin{displaymath}
	E^{sym}(\cM, \cN) = E^{inj}(\cM, \cN) + E^{inj}(\cN, \cM).
\end{displaymath}

We have that $E(\cM) = E^{inj}(\cM, \cM) = (1/2)E^{sym}(\cM, \cM)$.

%--------------------------------------------------------------------------------------
\subsection{Cluster categories and cluster characters}
In this section, we recall the definition of the (generalized) cluster category of a quiver with potential from \cite{Amiot08} and some results on cluster characters from \cite{Plamondon09}.

%................................
\subsubsection{Cluster categories}

Let $(Q,W)$ be any quiver with potential.  We first recall a construction of \cite{G06}.  Define a graded quiver $\overline{Q}$ from $Q$ in the following way:
\begin{itemize}
	\item the quiver $\overline{Q}$ has the same vertices as the quiver $Q$;
	\item the set of arrows of $Q$ is contained in that of $\overline{Q}$, and these arrows have degree $0$;
	\item for each arrow $a:i\rightarrow j$ of $Q$, add an arrow $a^{\ast}:j\rightarrow i$ of degree $-1$ in $\overline{Q}$;
	\item for each vertex $i$ of $Q$, add a loop $t_i:i\rightarrow i$ of degree $-2$ in $\overline{Q}$.
\end{itemize}

From the graded quiver $\overline{Q}$, we construct a differential graded algebra (dg algebra for short) $\Gamma = \Gamma_{Q,W}$ as follows.

As a graded algebra, $\Gamma$ is the completed path algebra of $\overline{Q}$.  In particular, its degree $\ell$ component is
\begin{displaymath}
	\Gamma^\ell = \prod_{\substack{w \textrm{ path in } \overline{Q} \\ deg(w) = \ell}} kw,
\end{displaymath}
for any integer $\ell$.

The differential $d$ of $\Gamma$ is the unique continuous and $k$-linear differential acting as follows on the arrows:
\begin{itemize}
	\item for any arrow $a$ of $Q$, $d(a) = 0$;
	\item for any arrow $a$ of $Q$, $d(a^\ast) = \partial_a W$;
	\item for any vertex $i$ of $Q$, $d(t_i) = e_i\big( \sum_{a\in Q_1}(aa^\ast - a^\ast a) \big)e_i$.
\end{itemize}

The dg algebra thus defined is the \emph{complete Ginzburg dg algebra}.  Notice that $\Ho^0\Gamma$ is canonically isomorphic to $J(Q,W)$.

Consider now the derived category $\cD\Gamma$ of $\Gamma$ (for background material on the derived category of a dg algebra, see, for example, \cite{K94} or \cite{KY09}).  Let $\perf\Gamma$ be the \emph{perfect derived category} of $\Gamma$, that is, the smallest triangulated full subcategory of $\cD\Gamma$ containing $\Gamma$ and closed under taking direct summands.  Denote by $\cD_{fd}\Gamma$ the full subcategory of $\cD\Gamma$ whose objects are those of $\cD\Gamma$ with finite-dimensional total homology.

It is proved in \cite[Theorem 2.17]{KY09} that $\cD_{fd}\Gamma$ is a triangulated subcategory of $\perf\Gamma$.

Following \cite[Definition 3.5]{Amiot08} (and \cite[Section 4]{KY09} in the case where $J(Q,W)$ is infinite-dimensional), we define the \emph{cluster category of $(Q,W)$} as the idempotent completion of the triangulated quotient $\perf \Gamma / \cD_{fd}\Gamma$.  We denote it by $\cC = \cC_{Q,W}$.

%................................
\subsubsection{Cluster characters}\label{sect::clusterchar}

Let $\cT$ be a triangulated category.  Let $T=\bigoplus_{i=1}^{n}T_i$ be a rigid object of $\cT$, where the $T_i$'s are indecomposable and pairwise non-isomorphic.  Assume that $\add T$ is a Krull--Schmidt subcategory.  Define the category $\pr{\cT}T$ as the full subcategory of $\cT$ whose objects are those $X$ for which there exists a triangle
\begin{displaymath}
	T_1^X \longrightarrow T_0^X \longrightarrow X \longrightarrow \Sigma T_1^X
\end{displaymath}  
with $T_0^X$ and $T_1^X$ in $\add T$.

In this situation, following \cite{DK08} and \cite{Palu08}, we define the \emph{index of $X$ with respect to $T$} as the element of the Grothendieck group $\K_0(\add T)$ given by
\begin{displaymath}
	\ind{T}X = [T_0^X] - [T_1^X].
\end{displaymath}
This is well-defined since we assumed $\pr{\cT}T$ to be Krull--Schmidt.

Our preferred example is the case when $\cT$ is the cluster category $\cC$ of a quiver with potential $(Q,W)$, and $T$ is the object $\Sigma^{-1}\Gamma$ (or any mutation of $\Sigma^{-1}\Gamma$ in the sense of \cite{KY09}, see section \ref{sect::derived}).  We assume now that we work in that case.

Consider the full subcategory $\cD$ of $\cC$ whose objects are those $M$ of $\pr{\cC}T \cap \pr{\cC}\Sigma T$ such that $\Hom{\cC}(T, M)$ is finite-dimensional.

Following \cite[Definition 3.10]{Plamondon09}, we define the \emph{cluster character with respect to $T$} as the map sending each (isomorphism class of) object $M$ of $\cD$ to the element (notations are explained below)
\begin{displaymath}
	X'_M = x^{\ind{T}\Sigma^{-1}M} \sum_{e \in \bN^n} \chi \Big( \Gr{e}\big( \Hom{\cC}(T,M)   \big)  \Big) x^{-\iota(e)}
\end{displaymath}
of $\bQ(x_1, \ldots, x_n)$.

Here, $\chi$ is the Euler-Poincar\'e characteristic, and for any $\End{\cC}(T)$-module $Z$, $\Gr{e}(Z)$ is the quiver Grassmannian of dimension vector $e$ of $Z$, that is, the projective variety whose points are the submodules of $Z$ with dimension vector $e$.  For any dimension vector $e$, $\iota(e)$ is the expression $\ind{T}Y + \ind{T}\Sigma^{-1}Y$ for any object $Y$ of $\cD$ such that $e$ is the dimension vector of $\Hom{\cC}(T, Y)$ (it was shown in \cite[Lemma 3.6]{Plamondon09} that $\iota(e)$ does not depend on the choice of such an object $Y$).  Finally, for any element $\underline{t} = \sum_{i=1}^n \lambda_i[T_i]$ of $\K_0(\add T)$, we denote by $x^{\underline{t}}$ the product $\prod_{i=1}^{n}x_i^{\lambda_i}$.

As proved in \cite[Proposition 3.6]{Plamondon09}, this cluster character satisfies the following \emph{multiplication formula}: for any objects $X$ and $Y$ of $\cD$ such that $\Hom{\cC}(X, \Sigma Y)$ is one-dimensional, if
\begin{displaymath}
	X\longrightarrow E \longrightarrow Y \longrightarrow \Sigma X \quad \textrm{and} \quad Y\longrightarrow E' \longrightarrow X \longrightarrow \Sigma Y
\end{displaymath}
are non-split triangles, then $E$ and $E'$ lie in $\cD$, and we have the equality 
\begin{displaymath}
	X'_X X'_Y = X'_E + X'_{E'}.
\end{displaymath}

%--------------------------------------------------------------------------------------
\subsection{Mutations as derived equivalences}\label{sect::derived}
Let $(Q,W)$ be a quiver with potential.  Assume that $Q$ has no loops, and that $i$ is a vertex of $Q$ not contained in a cycle of length $2$.  Let $(Q',W')$ be the mutated quiver with potential $\widetilde{\mu}_i(Q,W)$.

Let $\Gamma$ and $\Gamma'$ be the complete Ginzburg dg algebras associated with $(Q,W)$ and $(Q',W')$, respectively.

We recall here some results of \cite{KY09} on the mutation of $\Gamma$ in $\cD\Gamma$.  

Let $\Gamma_i^*$ be the cone in $\cD\Gamma$ of the morphism
\begin{displaymath}
	     \Gamma_i \longrightarrow \bigoplus_{\alpha}\Gamma_{t(\alpha)}
     \end{displaymath}
     whose components are given by left multiplication by $\alpha$.  Similarly, let $\Sigma\overline{\Gamma}_i^*$ be the cone of the morphism
  \begin{displaymath}
	      \bigoplus_{\beta}\Gamma_{s(\beta)} \longrightarrow \Gamma_i 
     \end{displaymath} 
     whose components are given by left multiplication by $\beta$.
     
Then it is proved in the discussion after \cite[Lemma 4.4]{KY09} that the morphism $\varphi_i:\Sigma\Gamma_i^* \longrightarrow \Sigma\overline{\Gamma}_i^*$   given in matrix form by
\begin{displaymath}
\left( \begin{array}{cc}
-\beta^{\ast} & -\partial_{\alpha\beta}W  \\
t_i & a^{\ast}  \\
\end{array} \right)
\end{displaymath}
becomes an isomorphism in $\cC$.

\begin{remark}\label{rema::W}
In particular, the composition of the morphisms
\begin{displaymath}
	\xymatrix{\bigoplus_{\alpha}\Gamma_{t(\alpha)} \ar[r] & \Gamma_i^* \ar[r]^{\Sigma^{-1}\varphi_i} & \overline{\Gamma}_i^* \ar[r] & \bigoplus_{\beta}\Gamma_{s(\beta)}
	}
\end{displaymath}
is given in matrix form by $\big( -\partial_{ab}W \big)$.
\end{remark}

The theorem below describes how the mutation of an indecomposable summand of $\Gamma$ can be interpreted as a derived equivalence.

\begin{theorem}[\cite{KY09}, Theorem 3.2]\label{theo::KY}
 \begin{enumerate}
   \item There exists a triangle equivalence $\widetilde{\mu}_i^+$ from $\cD(\Gamma')$ to $\cD(\Gamma)$ sending $\Gamma'_j$ to $\Gamma_j$ if \ $i \neq j$ and to the cone $\Gamma_i^*$ of the morphism 
     \begin{displaymath}
	     \Gamma_i \longrightarrow \bigoplus_{\alpha}\Gamma_{t(\alpha)}
     \end{displaymath}
     whose components are given by left multiplication by $\alpha$ if \ $i=j$.  The functor $\widetilde{\mu}_i^+$ restricts to triangle equivalences from $\perf \Gamma'$ to $\perf \Gamma$ and from $\cD_{fd} \Gamma'$ to $\cD_{fd} \Gamma$.
     
   \item Let $\Gamma_{red}$ and $\Gamma'_{red}$ be the complete Ginzburg dg algebra of the reduced part of $(Q,W)$ and $\tilde{\mu}_i(Q,W)$, respectively.  The functor $\widetilde{\mu}_i^+$ induces a triangle equivalence $(\mu_i^+)_{red} : \cD(\Gamma'_{red}) \longrightarrow \cD(\Gamma_{red})$ which restricts to triangle equivalences from $\perf \Gamma_{red}'$ to $\perf \Gamma_{red}$ and from $\cD_{fd} \Gamma'_{red}$ to $\cD_{fd} \Gamma_{red}$.
 \end{enumerate}
\end{theorem}
We will denote a quasi-inverse of $\widetilde{\mu}_i^+$ by the symbol $\widetilde{\mu}_i^-$.  Note that these equivalences induce equivalences on the level of cluster categories, which we will also denote by $\widetilde{\mu}_i^+$ and $\widetilde{\mu}_i^-$.

In Section \ref{sect::bijection}, we will need a concrete description of $\widetilde{\mu}_i^+$ and $\widetilde{\mu}_i^-$.  The functor $\widetilde{\mu}_i^+$ is the derived functor $? \otimes^{L}_{\Gamma'} T$, where $T$ is the $\Gamma'$-$\Gamma$-bimodule described below.  The functor $\widetilde{\mu}_i^-$ is then $\mathcal{H}om_{\Gamma}(T, ?)$.

As a right $\Gamma$-module, $T$ is a direct sum $\bigoplus_{j=1}^{n}T_j$, where $T_j$ is isomorphic to $e_j\Gamma$ if $i\neq j$ and $T_i$ is the cone of the morphism
\begin{displaymath}
	e_i\Gamma \ \longrightarrow \bigoplus_{\substack{\alpha\in Q_1 \\  s(\alpha)=i}}  e_{t(\alpha)}  \Gamma,
\end{displaymath}
whose components are given by left multiplication by $\alpha$.  Thus, as a graded module, $T_i$ is isomorphic to 
\begin{displaymath}
	P_{\Sigma i}\oplus \bigoplus_{\substack{\alpha\in Q_1 \\ s(\alpha)=i}}P_{\alpha},
\end{displaymath}
where $P_{\Sigma i}$ is a copy of $\Sigma (e_i\Gamma)$, and each $P_{\alpha}$ is a copy of $e_t(\alpha)\Gamma$.  We will denote by $e_{\Sigma i}$ the idempotent of $P_{\Sigma i}$ and by $e_{\alpha}$ the idempotent of $P_{\alpha}$.

The left $\Gamma'$-module structure of $T$ is described in terms of a homomorphism of dg algebras $\Gamma' \longrightarrow \mathcal{H}om_{\Gamma}(T,T)$, using the left $\mathcal{H}om_{\Gamma}(T,T)$-module structure of $T$.  We will need the description of the image of some elements of $\Gamma'$ under this homomorphism.  This description is given below.

For any vertex $j$ of $Q$, the element $e_j$ is sent to the identity of $T_j$.

Any arrow $\delta$ not incident with $i$ is sent to the map which is left multiplication by $\delta$.

For any arrow $\alpha$ of $Q$ such that $s(\alpha) = i$, the element $\alpha^{\star}$ is sent to the map $f_{\alpha^{\star}}:T_{t(\alpha)}\longrightarrow T_i$ defined by $f_{\alpha^{\star}}(a) = e_{\alpha} a$.

For any arrow $\beta$ of $Q$ such that $t(\beta) = i$, the element $\beta^{\star}$ is sent to the map $f_{\beta^{\star}}:T_i \longrightarrow T_{s(\beta)}$ defined by $f_{\beta^{\star}}(e_{\Sigma i}a_i + \sum_{s(\rho)=i}e_{\rho}a_{\rho})= -\beta^{\ast}a_i - \sum_{s(\rho)=i}(\partial_{\rho\beta}W)a_{\rho}$.

\section{Application to skew-symmetric cluster algebras}\label{sect::results}

%----------------------------------------------------------------------
\subsection{Rigid objects are determined by their index}\label{sect::rigid}
This section is the $\Hom{}$-infinite equivalent of \cite[Section 2]{DK08}.

Le $\cC$ be a triangulated category, and let $T=\bigoplus_{i=1}^{n}T_i$ be a rigid object of $\cC$, where the $T_i$'s are indecomposable and pairwise non-isomorphic.  Assume that $\pr{\cC}T$ is a Krull--Schmidt category, and that $B = \End{\cC}T$ is the completed Jacobian algebra $J(Q,W)$ of a quiver with potential $(Q,W)$.  An example of such a situation is the cluster category $\cC_{Q,W}$, with $T = \Sigma^{-1}\Gamma$.

The main result of this section is the following.

\begin{proposition}\label{prop::rigidindex}
With the above assumptions, if $X$ and $Y$ are rigid objects in $\pr{\cC}T$ such that $\ind{T}X = \ind{T}Y$, then $X$ and $Y$ are isomorphic.
\end{proposition}

The rest of the Section is devoted to the proof of the Proposition.

Let $X$ be an object of $\pr{\cC}T$, and let $\xymatrix{T_1^X \ar[r]^{f^X} & T_0^X \ar[r] & X \ar[r] & \Sigma T_1^X}$ be an $(\add T)$-presentation of $X$.  The group $\Aut{\cC}(T_1^X) \times \Aut{\cC}(T_0^X)$ acts  on the space $\Hom{\cC}(T_1^X, T_0^X)$, with action defined by $(g_1, g_0)f' = g_0f'(g_1)^{-1}$.  The orbit of $f^X$ under this action is the image of the map
\begin{eqnarray*}
	\Phi : \Aut{\cC}(T_1^X) \times \Aut{\cC}(T_0^X) & \longrightarrow & \Hom{\cC}(T_1^X, T_0^X) \\
	   (g_1, g_0) & \longmapsto & g_0f^X(g_1)^{-1}.
\end{eqnarray*}

Our strategy is to show that if $Y$ is another rigid object of $\pr{\cC}T$, then the orbits of $f^X$ and $f^Y$ must intersect (and thus coincide), proving that $X$ and $Y$ are isomorphic.

It was proved in \cite[Lemma 3.2]{Plamondon09} that the functor $F= \Hom{\cC}(T, ?)$ induces an equivalence of categories
\begin{displaymath}
	\pr{\cC}T/(\Sigma T) \longrightarrow \MOD B,
\end{displaymath}
where $\MOD B$ is the category of finitely presented right $B$-modules. Since $T$ is rigid, this implies that $F$ induces a fully faithful functor
\begin{displaymath}
 \add T \longrightarrow \MOD	B.
\end{displaymath} 
Thus we can often consider automorphisms and morphisms in the category $\MOD B$ instead of in $\cC$.

Now, let $\mathfrak{m}$ be the ideal of $J(Q,W)$ generated by the arrows of $Q$.  

The group $A = \Aut{B}(FT_1^X) \times \Aut{B}(FT_0^X)$ is the limit of the finite-dimensional affine algebraic groups 
\begin{displaymath}
 A_n = \Aut{B}(FT_1^X/(FT_1^X\mathfrak{m}^n)) \times \Aut{B}(FT_0^X/(FT_0^X\mathfrak{m}^n))
\end{displaymath}
with respect to the natural projection maps from $A_{n+1}$ to $A_{n}$, for $n \in \mathbb{N}$.

Similarly, the vector space $H = \Hom{B}(FT_1^X, FT_0^X)$ is the limit of the spaces
\begin{displaymath}
 H_n = \Hom{B}\Big(FT_1^X/(FT_1^X\mathfrak{m}^n), \ FT_0^X/(FT_0^X\mathfrak{m}^n)\Big)
\end{displaymath}
with respect to the natural projections.  All the $H_n$ are finite-dimensional spaces, and they are endowed with the Zariski topology.  The projection maps are then continuous, and $H$ is endowed with the limit topology.

Finally, for any integer $n$, we define a morphism $\Phi_n : A_n \rightarrow H_n$ which sends any element $(g_1, g_0)$ of $A_n$ to $g_0 f_n^X (g_1)^{-1}$, where $f_n^X$ is the image of $f^X$ in $H_n$ under the canonical projection.  Then the morphism $\Phi$ is the limit of the $\Phi_n$'s.

The situation is summarized in the following commuting diagram.

\begin{displaymath}
	\xymatrix{ A = \lim A_n\ar[d]^{\Phi} & \ldots \ar[r]& \ldots\ar[r] & A_3\ar[r]\ar[d]^{\Phi_3} & A_2\ar[r]\ar[d]^{\Phi_2} & A_1\ar[d]^{\Phi_1} \\
	           H = \lim H_n & \ldots \ar[r]& \ldots\ar[r] & H_3\ar[r] & H_2\ar[r] & H_1. 
	}
\end{displaymath}

The next step is the following : we will prove that the image of $\Phi$ is the limit of the images of the $\Phi_n$'s.

This will follow from the Lemma below.

\begin{lemma}\label{lemm::limits}
Let $(X_i)_{i\in \bN}$ be a family of topological spaces.  Let $(f_i:X_i\rightarrow X_{i-1})_{i\geq 1}$ be a family of continuous maps, and let $X = \Lim X_i$.  Let $(X'_i)_{i\in \bN}$ be another family of topological spaces, with continuous maps $(f'_i:X'_i\rightarrow X'_{i-1})_{i\geq 1}$, and let $X'= \Lim X'_i$.  Let $(u_i:X_i\rightarrow X'_i)$ be a familiy of continuous maps such that $f'_iu_i = u_{i-1}f_i$ for all $i\geq 1$, and let $u = \Lim u_i$.  Denote by $p_i: X\rightarrow X_i$ and $p'_i: X'\rightarrow X'_i$ the canonical projections.

For  integers $i<j$, denote by $f_{ij}$ (respectively $f'_{ij}$) the composition $f_jf_{j-1}\ldots f_{i+1}$ (respectively $f'_jf'_{j-1}\ldots f'_{i+1}$).  

Let $x'$ be an element of $X'$ with the property that for all $i\in \bN$, there exists $j \geq i$ such that for all $\ell \geq j$, $f_{i\ell}(u_{\ell}^{-1}(p'_{\ell}(x'))) = f_{ij}(u_{j}^{-1}(p'_{j}(x')))$.

Then $x'$ admits a preimage in $X$, that is, there exists $x\in X$ such that $u(x) = x'$.
\end{lemma}
\demo{ This is a consequence of the Mittag-Leffler theorem, see for instance \cite[Corollary II.5.2]{Bourbaki71}.
}

The above Lemma implies that the image of $\Phi$ is the limit of the images of the $\Phi_n$.  Indeed, the universal property of the limit gives an inclusion from the image of $\Phi$ to the limit of the images of the $\Phi_n$.  Let now $x'$ be in the image of $\Phi$, and let $x'_n$ be its projection in the image of $\Phi_n$.   The set $\Phi_n^{-1}(x_n)$ is a closed subset of $A_n$, and for any $m\geq n$, the image of $\Phi_m^{-1}(x_m)$ in $\Phi_n^{-1}(x_n)$ is closed.  Since $A_n$ has finite dimension as a variety, the sequences of images of the $\Phi_m^{-1}(x_m)$ in $\Phi_n^{-1}(x_n)$ eventually becomes constant.  Applying the above Lemma, we get that $x'$ has a preimage in $A$ by $\Phi$.  This proves that the image of $\Phi$ is the limit of the images of the $\Phi_n$.

We will now prove that the image of each $\Phi_n$ is open (and thus dense, since $H_n$ is irreducible).

To prove this, we pass to the level of Lie algebras.  To lighten notations, we let $E_n = \End{B}(FT_1^X/FT_1^X \mathfrak{m}) \times \End{B}(FT_0^X/FT_0^X \mathfrak{m})$ be the Lie algebra of $A_n$ for all positive integers $n$.  To prove that the image of $\Phi_n$ is open, it is sufficient to show that the map
\begin{eqnarray*}
	\Psi_n : E_n & \longrightarrow & H_n \\
	 (g_1, g_0) & \longmapsto & g_0f^X_n - f^X_ng_1
\end{eqnarray*}
is surjective.

The limit of the $E_n$'s is $E = \End{B}(FT_1^X) \times \End{B}(FT_0^X)$, and the limit of the $\Psi_n$'s is the map
\begin{eqnarray*}
	\Psi : E & \longrightarrow & H \\
	 (g_1, g_0) & \longmapsto & g_0f^X - f^Xg_1.
\end{eqnarray*}

The diagram below summarizes the situation.

\begin{displaymath}
	\xymatrix{ E = \lim E_n\ar[d]^{\Psi} & \ldots \ar[r]& \ldots\ar[r] &E_3\ar[r]\ar[d]^{\Psi_3} & E_2\ar[r]\ar[d]^{\Psi_2} & E_1\ar[d]^{\Psi_1} \\
	           H = \lim H_n & \ldots \ar[r]& \ldots\ar[r] & H_3\ar[r] & H_2\ar[r] & H_1. 
	}
\end{displaymath}

All the canonical projections are surjective.  

\begin{lemma}
The map $\Psi$ defined above is surjective. 
\end{lemma}
\demo{ This proof is contained in the proof of \cite[Lemma 2.1]{DK08} 
}

As a consequence, all the $\Psi_n$'s are surjective.  Hence the images of the $\Phi_n$'s are open.  

From this, we deduce that if $Y$ is another rigid object of $\pr{\cC}T$ with $(\add T)$-presentation $\xymatrix{T_0^X \ar[r]^{f^Y} & T_1^X\ar[r] & Y \ar[r] & \Sigma T_1^X}$, then $X$ and $Y$ are isomorphic.  Indeed, by the above reasonning, the orbit of $f^Y$ is the limit of the orbits of its projections in the $H_n$'s.  But these orbits are open, and so they intersect (and coincide) with the images of the $\Phi_n$ defined above.  Hence the orbit of $f^Y$ in $H$ is the limit of the images of the $\Phi_n$'s, and this is exactly the orbit of $f^X$.  Therefore $X$ and $Y$ are isomorphic.

The last step in proving Proposition \ref{prop::rigidindex} is to show that given $\ind{T}X$, we can ``deduce" $T_1^X$ and $T_0^X$.

An $(\add T)$-approximation $T_1^X \rightarrow T_0^X \rightarrow X \rightarrow \Sigma T_1^X$ is \emph{minimal} if one of the following conditions hold.
  \begin{itemize}
    \item The above triangle does not admit a direct summand of the form 
      \begin{displaymath}
        \xymatrix{R \ar[r]^{id_R}& R \ar[r] & 0\ar[r] & \Sigma R.}
      \end{displaymath}
    
    \item The morphism $f: T_0^X \rightarrow X$ in the presentation has the property that for any $g: T_0^X \rightarrow T_0^X$, the equality $fg = f$ implies that $g$ is an isomorphism.
  \end{itemize}

In fact, any of these two conditions implies the other.

\begin{lemma}
The above two conditions are equivalent if $\pr{\cC}T$ is Krull--Schmidt. 
\end{lemma}
\demo{ First suppose that the presentation has the form
\begin{displaymath}
	\xymatrix{ T'_1 \oplus R \ar[r]^{u\oplus 1_R} & T'_0 \oplus R \ar[r]^{\phantom{xx}(f', 0)} & X \ar[r] & \Sigma T_1^X,
}
\end{displaymath}
where $f = (f', 0)$ in matrix form.

Then the endomorphism $g$ of $T'_0 \oplus R$ given by $g = 1_{T'_0} \oplus 0$ is not an isomorphism, and $fg = f$.

Now suppose that the presentation admits no direct summand of the form
\begin{displaymath}
   \xymatrix{R \ar[r]^{id_R}& R \ar[r] & 0\ar[r] & \Sigma R.}
\end{displaymath}
Using the Krull--Schmidt property of $\pr{\cC}T$, we can decompose bot $T_0^X$ and $T_1^X$ as a finite direct sum of objects with local endomorphism rings.  In that case, the morphism $f$ written in matrix form (in any basis) has no non-zero entries.

Let $g$ be an endomorphism of $T_0^X$ such that $fg = f$.  Then $f(1_{T_0^X} - g) = 0$.  Consider the morphism $(1_{T_0^X} - g)$ written in matrix form.  If one of its entries is an isomorphism, then by a change of basis we can write $(1_{T_0^X} - g)$ as the matrix
\begin{displaymath}
 \left(\begin{array}{c|c}
 * & 0 \\
 \hline
 0 & \phi
 \end{array}\right),
 \end{displaymath}
 where $\phi$ is an isomorphism.  In that case, it is impossible that $f(1_{T_0^X} - g)=0$, since $f$ has no non-zero entries. This implies that none of the entries of the matrix of $(1_{T_0^X} - g)$ is invertible.  Therefore the diagonal entries of $g$ are invertible (since for any element $x$ of a local ring, if $(1-x)$ is not invertible, then $x$ is), while the other entries are not, and $g$ is an isomorphism.
}

\begin{lemma}\label{lemm::repetition}
If $X$ is rigid and $T_1^X \stackrel{\alpha}{\rightarrow} T_0^X \rightarrow X \stackrel{\gamma}{\rightarrow} \Sigma T_1^X$ is a minimal $(\add T)$-presentation, then $T_1^X$ and $T_0^X$ have no direct summand in common.
\end{lemma}
\demo{The first proof of \cite[Proposition 2.2]{DK08} works in this setting. We include here a similar argument for the convenience of the reader.

Suppose that $T_i$ is a direct factor of $T_0^X$.  Let us prove that it is not a direct factor of $T_1^X$.

Applying $F = \Hom{\cC}(T, ?)$ to the triangle above, we get a minimal projective presentation of $FX$.  This yields an exact sequence
\begin{displaymath}
	(FX, S_i) \longrightarrow (FT_0^X, S_i) \stackrel{F\alpha^*}{\longrightarrow} (FT_1^X, S_i),
\end{displaymath}
Where $S_i$ is the simple at the vertex $i$. Since the presentation is minimal, $F\alpha^*$ vanishes, and there exists a non-zero morphism $f:FX\rightarrow S_i$.  In particular, $f$ is surjective.

Let $g:FT_1^X \longrightarrow S_i$ be a morphism. Since $FT_1^X$ is projective, there exists a morphism $h:FT_1^X \longrightarrow FX$ such that $fh = g$.

Lift $S_i$ to an object $\Sigma T_i^*$ of $\cC$, and lift $f$, $g$, and $h$ to morphisms $\overline{f}:X\rightarrow \Sigma T_i^*$, $\overline{g}:T_1^X\rightarrow \Sigma T_i^*$ and $\overline{h}:T_1^X\rightarrow X$ of $\cC$.  Using \cite[Lemma 3.2 (1)]{Plamondon09}, we get that $\overline{f}\overline{h} = \overline{g}$.
\begin{displaymath}
  \xymatrix{ \Sigma^{-1}X \ar[r]^{\Sigma^{-1}\gamma} & T_1^X \ar[r]^{\alpha}\ar[d]_{\overline{g}}\ar[dl]_{\overline{h}} & T_0^X\ar@{.>}[dl]_{\sigma} \\
             X \ar[r]^{\overline{f}} & \Sigma T_i^*
  }	
\end{displaymath}
Since $X$ is rigid, $\overline{h}\Sigma^{-1}\gamma$ vanishes, and thus so does $\overline{g}\Sigma^{-1}\gamma$.  Then there exists a morphism $\sigma: T_0^X \rightarrow \Sigma T_i^*$ such that $\sigma\alpha = \overline{g}$.  But since $F\alpha^* = 0$, we get that $g = (F\sigma)(F\alpha)$ vanishes.

We have thus shown that there are no non-zero morphisms from $FT_1^X$ to $S_i$.  Therefore $T_i$ is not a direct factor of $T_1^X$.
}

By the above Lemma, the knowledge of $\ind{T}X$ is sufficient to deduce the isomorphism classes of $T_1^X$ and $T_0^X$ in any minimal $(\add)$-presentation of $X$.  Therefore, if $Y$ is another rigid object of $\pr{\cC}T$ with $\ind{T}X = \ind{T}Y$, all of the above reasonning implies that $X$ and $Y$ are isomorphic.  This finishes the proof of Proposition \ref{prop::rigidindex}.  

%------------------------------------------------------------
\subsection{Index and $\bg$-vectors}\label{sect::index g-vectors}

It was proved in \cite[Proposition 6.2]{FK09} that, inside a certain $\Hom{}$-finite cluster category $\cC$, the index of an object $M$ with respect to a cluster-tilting object $T$ gives the $\bg$-vector of $X'_M$ with respect to the associated cluster.  The authors then used this result to prove conjectures of \cite{FZ07} in this case.

In this section, we will prove a similar result, dropping the assumption of $\Hom{}$-finiteness.

Let $(Q,F)$ be a finite ice quiver, where $Q$ has no oriented cycles of length $\leq 2$.  Suppose that the associated matrix $B$ has full rank $r$.  Denote by $\cA$ the associated cluster algebra.  Let $W$ be a potential on $Q$, and let $\cC = \cC_{Q,W}$ be the associated cluster category.  Denote by $\cD$ the full subcategory of $\pr{\cC}\Gamma \cap \pr{\cC} \Sigma^{-1}\Gamma$ whose objects are those $X$ such that $\Hom{\cC}(\Sigma^{-1} \Gamma, X)$ is finite-dimensional.

Following \cite{FK09}, let $\cU$ be the full subcategory of $\cD$ defined by 
\begin{displaymath}
	\cU = \{X \in \cD \ \big| \ \Hom{\cC}(\Sigma^{-1}\Gamma_j, X)=0 \ \textrm{for } r+1 \leq j \leq n \}.
\end{displaymath}
Note that $\cU$ is invariant under iterated mutation of $\Gamma$ at vertices $1, 2, \ldots, r$.

Let $T=\bigoplus_{j=1}^n T_j = \bigoplus_{j=1}^r T_j \oplus \bigoplus_{j=r+1}^n \Gamma_j$ be a rigid object of $\cD$ reachable from $\Gamma$ by mutation at an admissible sequence of vertices of $Q$ not in $F$, and let $G$ be the functor $\Hom{\cC}(\Sigma^{-1}T, ?)$ from $\cC$ to the category of $\End{\cC}(T)$-modules.  Let $X'_?$ be the associated cluster character, defined by
\begin{displaymath}
	X'_M = x^{\ind{T}M}\sum_e \Big( \chi\big( \Gr{e}(GM) \big) \Big) x^{-\iota(e)},
\end{displaymath}
where $\iota(e)$ is the vector $\ind{T}Y + \ind{T}\Sigma Y$ for any $Y$ such that the dimension vector of $FY$ is $e$ (it was proved in \cite[Lemma 3.6]{Plamondon09} that this vector is independent of the choice of such a $Y$, see also \cite{Palu08}).

 Since we only allow mutations at vertices not in $F$, the Gabriel quiver of $T$ can be thought of as an ice quiver $(Q^T, F)$ with same set of frozen vertices as $(Q,T)$.  Let $B^T = (b_{j\ell}^T)$ be the matrix associated to $(Q^T, F)$.  According to \cite[Lemma 1.2]{GSV03} and \cite[Lemma 3.2]{BFZ05} , $B^T$ is of full rank $r$ if $B$ is.
 
 Suppose now that $M$ is an object of $\cU$.  Let us prove that $X'_M$ then admits a $\bg$-vector, that is, $X'_M$ is in the set $\cM$ defined in Section \ref{sect::g-vectors}.  In order to do this, let us compute $-\iota(\delta_j)$, where $\delta_j$ is the vector whose $j$-th coordinate is $1$ and all others are $0$, for $j=1,2,\ldots, r$.
 
 Let $T_j^{\ast}$ be an indecomposable object of $\cD$ such that $GT_j^{\ast}$ is the simple $\End{\cC}(T)$-module at the vertex $j$. It follows from the derived equivalence in \cite[Theorem 3.2]{KY09} that we have triangles
 \begin{displaymath}
   T_j \rightarrow \bigoplus_{\substack{\alpha \in Q^T_1 \\ s(\alpha)=j}}T_{t(\alpha)} \rightarrow T_j^{\ast} \rightarrow \Sigma T_j	\quad \textrm{and } \quad T_j^{\ast} \rightarrow \bigoplus_{\substack{\alpha \in Q^T_1 \\ t(\alpha)=j}}T_{s(\alpha)} \rightarrow T_j \rightarrow \Sigma T_j^{\ast}.
 \end{displaymath}
 
 We deduce from those triangles that for any $0 \leq \ell \leq n$, the $\ell$-th entry of $-\iota(\delta_j)$ is the number of arrows in $Q^T$ from $\ell$ to $j$ minus the number of arrows from $j$ to $\ell$.  This number is $b_{\ell j}^T$.  Thus, with the notations of Section \ref{sect::g-vectors}, we have that $x^{-\iota(\delta_j)} = \prod_{\ell =1}^n x_{\ell}^{b_{\ell j}} = \hat{y}_j$.
 
 Therefore, since $\iota$ is additive, for $M$ in $\cU$, we have the equality
 \begin{displaymath}
   X'_M	= x^{\ind{T}M}\sum_e \Big( \chi\big( \Gr{e}(GM) \big) \Big)\prod_{j=1}^r \hat{y}_j^{e_j}
 \end{displaymath}
 (notice that if $M$ is in $\cU$, then $\Gr{e}(GM)$ is empty for all vectors $e$ such that one of $e_{r+1}, \ldots, e_n$ is non-zero).
 
 Moreover, the rational function $R(u_1, \ldots, u_r) = \sum_e \Big( \chi\big( \Gr{e}(GM) \big) \Big)\prod_{j=1}^r u_j^{e_j}$ is in fact a polynomial with constant coefficient $1$, and is thus primitive.  
 
 We have proved the following result.

\begin{proposition}\label{prop::index}
  Any object $M$ of $\cU$ is such that $X'_M$ admits a $\bg$-vector.  This $\bg$-vector $(g_1, \ldots, g_r)$ is given by $g_j = [\ind{T}M: T_j]$, for $1\leq j \leq r$.
\end{proposition}

These considerations allow us to prove the following Theorem, whose parts (1), (3) and (4) were first shown in the same generality in \cite{DWZ09} using decorated representations, and then in \cite{Nag09} using Donaldson--Thomas theory.

We say that a collection of vectors of $\bZ^r$ are \emph{sign-coherent} if the $i$-th coordinates of all the vectors of the collection are either all non-positive of all non-negative.

\begin{theorem}\label{theo::gvectors}
  Let $(Q,F)$ be any ice quiver without oriented cycles of length $\leq 2$, and let $\cA$ be the associated cluster algebra.  Suppose that the matrix $B$ associated with $(Q,F)$ is of full rank $r$.
\begin{enumerate}
	\item Conjecture 6.13 of \cite{FZ07} holds for $\cA$, that is, the $\bg$-vectors of the cluster variables of any given cluster are sign-coherent.
	\item Conjecture 7.2 of \cite{FZ07} holds for $\cA$, that is, the cluster monomials are linearly independent over $\bZ\bP$, where $\bP$ is the tropical semifield in the variables $x_{r+1}, \ldots, x_n$.
	\item Conjecture 7.10 of \cite{FZ07} holds for $\cA$, that is, different cluster monomials have different $\bg$-vectors, and the $\bg$-vectors of the cluster variables of any cluster form a $\bZ$-basis of $\bZ^r$.
	\item Conjecture 7.12 of \cite{FZ07} holds for $\cA$, that is, if $\bg = (g_1, \ldots, g_r)$ and $\bg' = (g'_1, \ldots, g'_r)$ are the $\bg$-vectors of one cluster monomial with respect to two clusters $t$ and $t'$ related by one mutation at the vertex $i$, then we have
	\begin{displaymath}
	  g'_j = \left\{ \begin{array}{ll}
                      -g_i & \textrm{if } j=i\\
                      g_j + [b_{ji}]_{+}g_{i} - b_{ji}\min(g_{i}, 0)  & \textrm{if } j\neq i
                   \end{array} \right.
  \end{displaymath}
  where $B = (b_{j\ell})$ is the matrix associated with the seed $t$, and we set $[x]_{+} = \max(x,0)$ for any real number $x$.
\end{enumerate}
\end{theorem}
\demo{ Choose a non-degenerate potential $W$ on $Q$, and let $\cC = \cC_{Q,W}$ be the associated cluster category. Let $X'_?$ be the cluster character associated with $\Gamma$. 

We first prove Conjecture 6.13.  We reproduce the arguments of \cite[Section 2.4]{DK08}.  To any cluster $t$ of $\cA$, we associate (using \cite[Theorem 4.1]{Plamondon09}) a reachable rigid object $T$ of $\cU$, obtained by mutating at vertices not in $F$.  Write $T$ as the direct sum of the indecomposable objects $T_1, \ldots, T_n$.    Then, for $1\leq j \leq r$, we have that $X'_{T_j}$ is a cluster variable lying in the cluster $t$.  By Proposition \ref{prop::index}, its $\bg$-vector $(g_1^j, \ldots, g_r^j)$ is given by $g^j_\ell = [\ind{\Gamma}T_j : \Gamma_\ell]$.  Now, by Lemma \ref{lemm::repetition}, any minimal $\add \Gamma$-presentation of $T$ 
\begin{displaymath}
	R_1 \longrightarrow R_0 \longrightarrow T \longrightarrow \Sigma R_1
\end{displaymath}
is such that $R_0$ and $R_1$ have no direct factor in common.  But this triangle is a direct sum of minimal presentations of $T_1, \ldots, T_n$.  Therefore the indices of these objects must be sign-coherent.  This proves Conjecture 6.13.

Next, we prove Conjecture 7.2.  We prove it in the same way as in \cite[Corollary 4.4 (b) and Theorem 6.3 (c)]{FK09}.  Using \cite[Theorem 4.1]{Plamondon09}, we associate to any finite collection of clusters $(t_j)_{j\in J}$ of $\cA$ a family of reachable rigid objects $(T^j)_{j\in J}$ of $\cU$, obtained by mutating at vertices not in $F$ (for the moment we do not know if this assignment is uninque).  Let $(M_j)_{j\in J}$ be a family of pairwise non-isomorphic objects, where each $M_j$ lies in $\add T^j$ (in particular, these objects are rigid).  Any $\bZ\bP$-linear combination of cluster monomials can be written as a $\bZ$-linear combination of some $X'_{M_j}$'s, where the $M_j$'s are as above.  Thus it is sufficient to show that the $M_j$'s are linearly indepenant over $\bZ$.

The key idea is to assign a degree to each $x_j$ in such a way that each $\hat{y}_j$ is of degree $1$.  Such an assignment is obtained by putting $\deg(x_j) = k_j$, where the $k_j$'s are rational numbers such that
\begin{displaymath}
	(k_1, \ldots, k_n) B = (1, \ldots, 1).
\end{displaymath} 
This equation admits a solution, since the rank of $B$ ir $r$.  Thus the term of minimal degree in $X'_M$ is $x^{\ind{\Gamma}M}$, for any $M$ in $\cU$.

Now let $(c_j)_{j\in J}$ be a family of real numbers such that $\sum_{j\in J}c_jX'_{M_j} = 0$.  The term of minimal degree of this polynomial has the form $\sum_{\ell\in L}c_{\ell}x^{\ind{\Gamma}M_{\ell}}$ for some non-empty subset $L$ of $J$, and this term must vanish.  But according to Proposition \ref{prop::rigidindex}, the indices of the $M_{\ell}$'s are pairwise distinct.  Thus $c_{\ell}$ is zero for any $\ell \in L$.  Repeating this argument, we get that $c_j$ is zero for any $j \in J$.  This proves the linear independance of cluster monomials.

%.................

The proof of Conjecture 7.10 goes as follows.  Let $\{ w_1, \ldots, w_r \}$ be a cluster of $\cA$, and let $w_1^{a_1} \ldots w_r^{a_r}$ be a cluster monomial.  Let $T = \bigoplus_{j=1}^{r}T_j \oplus \bigoplus_{j=r+1}^{n}\Gamma_n$ be the rigid object of $\cC$ associated with that cluster.  Then the cluster character
\begin{displaymath}
	X'_M = x^{\ind{\Gamma}M}\sum_e \Big( \chi\big( \Gr{e}(\Hom{\cC}(\Sigma^{-1}\Gamma, M) \big) \Big) x^{-\iota(e)}
\end{displaymath}
sends the object $\bigoplus_{j=1}^{r}T_j^{a_j}$ to the cluster monomial $w_1^{a_1} \ldots w_r^{a_r}$.  The $\bg$-vector of this cluster monomial is the index of $\bigoplus_{j=1}^{r}T_j^{a_j}$ by Proposition \ref{prop::index}, and by Proposition \ref{prop::rigidindex}, this object is completely determined by its index.  Therefore two different cluster monomials, being associated with different rigid objects of $\cC$, have different $\bg$-vectors.  

Let us now prove that the $\bg$-vectors of $w_1, \ldots, w_r$ form a basis of $\bZ^r$.  For any object $M$ of $\cD$, denote by $(\ind{\Gamma}M)_0$ the vector containing the first $r$ components of $\ind{\Gamma}M$.  In view of Proposition \ref{prop::index}, it is sufficient to prove that the vectors $(\ind{\Gamma}T_1)_0, \ldots, (\ind{\Gamma}T_r)_0$ form a basis of $\bZ^r$.  

We prove this by induction.  The statement is trivially true for $\Gamma$.  Now suppose it is true for some reachable object $T$ as above.  Let $1 \leq \ell \leq r$ be a vertex of $Q$, and let $T' = \mu_{\ell}(T)$.  We can write $T' = \bigoplus_{j=1}^{n} T'_j$, where $T'_j = T_j$ if $j\neq \ell$, and there are triangles
\begin{displaymath}
	T_{\ell} \longrightarrow \bigoplus_{\substack{\alpha\in Q^T_1 \\ s(\alpha) = \ell}} T_{t(\alpha)} \longrightarrow T'_{\ell} \longrightarrow \Sigma T_{\ell} \quad \textrm{and} \quad T'_{\ell} \longrightarrow \bigoplus_{\substack{\alpha\in Q^T_1 \\ t(\alpha) = \ell}} T_{s(\alpha)} \longrightarrow T_{\ell} \longrightarrow \Sigma T'_{\ell}
\end{displaymath}
thanks to \cite{KY09}.  Moreover, the space $\Hom{\cC}(T'_{\ell}, \Sigma T_{\ell})$ is one-dimensional; by applying \cite[Lemma 3.8]{Plamondon09} (with the $T$ of the Lemma being equal to our $\Sigma^{-1}\Gamma$), we get an isomorphism
\begin{displaymath}
	(\Gamma)(T'_{\ell}, \Sigma T_{\ell}) \longrightarrow D\Hom{\cC}(T_{\ell}, \Sigma T'_{\ell}) / (\Gamma).
\end{displaymath}
Therefore one of the two morphisms $T'_{\ell} \rightarrow \Sigma T_{\ell}$ and $T_{\ell} \rightarrow \Sigma T'_{\ell}$ in the triangles above is in $(\Gamma)$.  Depending on which one is in $(\Gamma)$, and applying \cite[Lemma 3.4(2)]{Plamondon09}, we get that either
\begin{displaymath}
	  \ind{\Gamma}T'_j = \left\{ \begin{array}{ll}
                      \ind{\Gamma}T_j & \textrm{if } j\neq\ell\\
                      -\ind{\Gamma}T_{\ell}  + \sum_{\substack{\alpha\in Q^T_1 \\ s(\alpha) = \ell}} \ind{\Gamma}T_{t(\alpha)}& \textrm{if } j=\ell.
                   \end{array} \right.
\end{displaymath}
or
\begin{displaymath}
	  \ind{\Gamma}T'_j = \left\{ \begin{array}{ll}
                      \ind{\Gamma}T_j & \textrm{if } j\neq\ell\\
                      -\ind{\Gamma}T_{\ell}  + \sum_{\substack{\alpha\in Q^T_1 \\ t(\alpha) = \ell}} \ind{\Gamma}T_{s(\alpha)}& \textrm{if } j=\ell.
                   \end{array} \right.
\end{displaymath}
holds.  Therefore the $(\ind{\Gamma}T'_j)_0$'s still form a basis of $\bZ^r$.  Conjecture 7.10 is proved.

%....................

Finally, let us now prove Conjecture 7.12.  Let $T$ and $T'$ be reachable rigid objects related by a mutation at vertex $\ell$, as above.  Then we have two triangles
\begin{displaymath}
	T_{\ell} \longrightarrow E \longrightarrow T'_{\ell} \longrightarrow \Sigma T_{\ell} \quad \textrm{and} \quad T'_{\ell} \longrightarrow E' \longrightarrow T_{\ell} \longrightarrow \Sigma T'_{\ell},
\end{displaymath}
where $E = \bigoplus_{\substack{\alpha\in Q^T_1 \\ s(\alpha) = \ell}} T_{t(\alpha)}$ and $E' = \bigoplus_{\substack{\alpha\in Q^T_1 \\ t(\alpha) = \ell}} T_{s(\alpha)}$.  Moreover, the dimension of the space $\Hom{\cC}(T, \Sigma T')$ is one.  Thus we can apply \cite[Proposition 2.7]{Plamondon09}. 

Let $M$ be a rigid object in $\pr{\cC}T$, and let $T_1^M \rightarrow T_0^M \rightarrow M \rightarrow \Sigma T_1^M$ be a minimal $(\add T)$-presentation.  Then, by \cite[Proposition 2.7]{Plamondon09}, $M$ is in $\pr{\cC}T'$.  Moreover, if $T_0^M = \overline{T}_0^M \oplus T_{\ell}^a$ and $T_1^M = \overline{T}_1^M \oplus T_{\ell}^b$, where $T_{\ell}$ is not a direct summand of $\overline{T}_0^M \oplus \overline{T}_1^M$, then the end of the proof of that Proposition gives us a triangle
\begin{displaymath}
	E^b \oplus (T'_{\ell})^{a-c} \oplus \overline{T}_1^{M} \longrightarrow (T'_{\ell})^{b-c} \oplus \overline{T}_0^{M} \oplus (E')^a \longrightarrow M \longrightarrow \Sigma(E^b \oplus (T'_{\ell})^{a-c} \oplus \overline{T}_1^{M}),
\end{displaymath}
for some integer $c$.  Notice that $[\ind{T}M : T_{\ell}] = (a-b)$, and that since $T_0^M$ and $T_0^M$ have no direct factor in common by Lemma \ref{lemm::repetition}, one of $a$ and $b$ must vanish.   Notice further that $b = -\min([\ind{T}M : T_{\ell}], 0)$.   Thus
\begin{displaymath}
	  [\ind{T'}M : T'_j] = \left\{ \begin{array}{ll}
                       -[\ind{T}M : T_{\ell}] \qquad (\textrm{if } j = \ell) \\
                        \phantom{}[\ind{T}M : T_{j}] + [\ind{T}M : T_{\ell}] [b^T_{j\ell}]_{+} - b^T_{j\ell}\min([\ind{T}M : T_{\ell}], 0) \\
                         \qquad (\textrm{if } j\neq\ell).
                   \end{array} \right.
\end{displaymath}
This proves the desired result on $\bg$-vectors. 
}

\begin{remark}\label{rema::indices}
Using the notations of the end of the proof of Theorem \ref{theo::gvectors}, we get that, if $M$ is an object of $\cD$ which is not necessarily rigid, then
\begin{displaymath}
	  [\ind{T'}M : T'_j] = \left\{ \begin{array}{ll}
                       -[\ind{T}M : T_{\ell}] \qquad (\textrm{if } j = \ell) \\
                        \phantom{}[\ind{T}M : T_{j}] + a[b_{j\ell}]_+ -b[-b_{j\ell}]_+ \qquad (\textrm{if } j\neq\ell).
                   \end{array} \right.
\end{displaymath}
Moreover, if the presentation $T_1^M \rightarrow T_0^M \rightarrow M \rightarrow \Sigma T_1^M$ is minimal, then the integer $c$ vanishes.  Indeed, in the proof of \cite[Proposition 2.7]{Plamondon09}, $c$ (or $r$ in \cite{Plamondon09}) is defined by means of the composition
\begin{displaymath}
	T_1^M \longrightarrow \overline{T}_0^M \oplus T_\ell^a \longrightarrow \Sigma(T'_\ell)^a.
\end{displaymath}
The minimality of the presentation implies that this composition vanishes, and thus that $c=0$.
\end{remark}

Using Theorem \ref{theo::gvectors}, we get a refinement of \cite[Theorem 4.1]{Plamondon09}.

\begin{corollary}
  The cluster character $X'_?$ associated with $\Gamma$ induces a bijection between the set of isomorphism classes of indecomposable reachable rigid objects of $\cC$ and the set of cluster variables of $\cA$.
\end{corollary}
\demo{  It was proved in \cite[Theorem 4.1]{Plamondon09} that we have a surjection.  We deduce from Theorem \ref{theo::gvectors} that different indecomposable reachable rigid objects are sent to different cluster variables.  Indeed, different such objects are sent to elements in $\cA$ which are linearly independent, and thus different.
}

We also get that the mutation of rigid objects governs the mutation of tropical $Y$-variables, as shown in \cite[Corollary 6.9]{K10} in the $\Hom{}$-finite case.

\begin{corollary}\label{coro::tropical}
Let $(Q,W)$ be a quiver with potential, and let $\cC$ be the associated cluster category.  Let $\underline{i} = (i_1, \ldots, i_m)$ be an admissible sequence of vertices, and let $T'$ be the object $\mu_{\underline{i}}(\Gamma)$.  Let $(Q, \by)$ be a $Y$-seed, with $\by = (y_1, \ldots, y_n)$.

Then $\mu_{\underline{i}}(Q, \by)$ is given by $(\mu_{\underline{i}}(Q), \by')$, where
\begin{displaymath}
	y_j' = \prod_{s=1}^{n} y_{s}^{-[\ind{\Sigma^{-1}T'}\Gamma_{s} : \Sigma^{-1}T'_j]}.
\end{displaymath}
\end{corollary}
\demo{ The result is proved by induction on $m$.  It is trivially true for $m=0$, that is, for empty sequences of mutations.  Suppose it is true for any sequence of $m$ mutations.

Let $\underline{i}' = (i_1, \ldots, i_m, \ell)$ be an admissible sequence of $m+1$ mutations.  Let $T'' = \mu_{\underline{i}'}(\Gamma)$ and $(\mu_{\underline{i}'}(Q), \by'') = \mu_{\underline{i}'}(Q, \by)$.

Using the mutation rule for $Y$-seeds (see section \ref{sect::tropical}) and the induction hypothesis, we get that
\begin{displaymath}
	y''_\ell = \prod_{s=1}^{n} y_{s}^{[\ind{\Sigma^{-1}T'}\Gamma_{s} : \Sigma^{-1}T'_j]}
\end{displaymath}
and that, for any vertex $j$ different from $\ell$, 
\begin{displaymath}
	y''_j = 
 \prod_{s=1}^{n} y_{s}^{-[\ind{\Sigma^{-1}T'}\Gamma_{s} : \Sigma^{-1}T'_j] - [\ind{\Sigma^{-1}T'}\Gamma_{s} : \Sigma^{-1}T'_\ell][b^{T'}_{\ell j}]_{+} - b^{T'}_{\ell j} \min(-[\ind{\Sigma^{-1}T'}\Gamma_{s} : \Sigma^{-1}T'_\ell] , 0)} 
\end{displaymath}

Now, recall from the end of the proof of Theorem \ref{theo::gvectors} that for any object $M$ of $\pr{\cC}T'$, we have an $(\add T'')$-presentation
\begin{displaymath}
	E^b \oplus (T''_{\ell})^{a-c} \oplus \overline{T'}_1^{M} \longrightarrow (T''_{\ell})^{b-c} \oplus \overline{T'}_0^{M} \oplus (E')^a \longrightarrow M \longrightarrow \Sigma(E^b \oplus (T''_{\ell})^{a-c} \oplus \overline{T'}_1^{M}),
\end{displaymath}
and that $[\ind{T'}M : T'_\ell] = (a-b)$.  Notice also that $a = -\min([-\ind{T'}M:T'_\ell],0)$.  Thus
\begin{displaymath}
	  [\ind{T''}M : T''_j] = \left\{ \begin{array}{ll}
                       -[\ind{T'}M : T'_{\ell}] \qquad (\textrm{if } j = \ell) \\
                        \phantom{}[\ind{T'}M : T'_{j}] + [\ind{T'}M : T'_{\ell}] [b^{T'}_{\ell j}]_{+} + \\
                         \qquad + b^{T'}_{\ell j}\min(-[\ind{T'}M : T'_{\ell}], 0)  \qquad (\textrm{if } j\neq\ell).
                   \end{array} \right.
\end{displaymath}
Replacing $M$ by $\Sigma \Gamma_s$, and using the above computation of $y''_j$, we get exactly the desired equality. 
}

\begin{remark}
The opposite category $\cC^{op}$ is triangulated with suspension functor $\Sigma_{op} = \Sigma^{-1}$.  If $T$ is a rigid object of $\cC$, then it is rigid in $\cC^{op}$, and any object $X$ admitting an $(\add \Sigma^{-1}T)$-presentation in $\cC$ admits an $(\add T)$-presentation in $\cC^{op}$.  If we denote by $\ind{T}^{op}X$ the index of $X$ with respect to $T$ in $\cC^{op}$, then we have the equality $\ind{T}^{op}X = -\ind{\Sigma^{-1}T}X$.  Thus the equality of Corollary \ref{coro::tropical} can be written as
\begin{displaymath}
	y_j' = \prod_{s=1}^{n} y_{s}^{[\ind{T'}^{op}\Gamma_{s} : T'_j]}.
\end{displaymath}
This corresponds to the notation and point of view adopted in \cite[Corollary 6.9]{K10}.
\end{remark}

%--------------------------------------------------------------------
\subsection{Cluster characters and $F$-polynomials}\label{sect::clustercharF}

Let $\cA$ be a cluster algebra with principal coefficients at a seed $\big( (Q,F), \bx \big)$.  In particular, $n=2r$, and the matrix $B$ associated with $(Q,F)$ has full rank $r$.

Let $W$ be a potential on $Q$, and let $\cC = \cC_{Q,W}$ be the cluster category associated with $(Q,W)$.  Let $T$ be a rigid object of $\cC$ reachable from $\Gamma$ by mutation at an admissible sequence of vertices $(i_1, \ldots, i_s)$ not in $F$.  Write $T$ as $\bigoplus_{j=1}^{2r}T_j$, where $T_{\ell} = \Gamma_{\ell}$ for $r < \ell \leq 2r$.

For any vertex $j$ not in $F$, $X'_{T_j}$ is a cluster variable in $\cA$.  Specializing at $x_1 = \ldots = x_r = 1$, we obtain the corresponding $F$-polynomial (see Section \ref{sect::F-polynomials}), which we will denote by $F_{T_j}$.

We thus have the equality
\begin{displaymath}
	F_{T_j} = \prod_{i=r+1}^{2r} x_i^{[\ind{\Gamma}T_j : \Gamma_i]} \sum_{e}\chi\Big( \Gr{e} \big( \Hom{\cC}(\Sigma^{-1}\Gamma, T_j) \big) \Big) \prod_{i=r+1}^{2r} x_i^{-\iota(e)_i},
\end{displaymath}
where $\iota(e)$ was defined in section \ref{sect::clusterchar} and $\iota(e)_i$ is the $i$-th component of $\iota(e)$.

\begin{remark}
The element $X'_{T_j}$ of $\cA$ is the $j$-th cluster variable of the cluster obtained from the initial cluster at the sequence of vertices $(i_1, \ldots, i_s)$ by \cite[Theorem 4.1]{Plamondon09}.  Therefore, the polynomial $F_{T_j}$ is the corresponding $F$-polynomial.
\end{remark}

It follows from our computation in Section \ref{sect::index g-vectors} that for $r<i\leq 2r$, there is an equality $-\iota(e)_i = \sum_{j=1}^r e_j b_{ij}$, and since our cluster algebra has principal coefficients, this number is $e_{i-r}$.  Thus we get the equality
\begin{displaymath}
	F_{T_j} = \prod_{i=r+1}^{2r} x_i^{[\ind{\Gamma}T_j : \Gamma_i]} \sum_{e}\chi\Big( \Gr{e} \big( \Hom{\cC}(\Sigma^{-1}\Gamma, T_j) \big) \Big) \prod_{i=r+1}^{2r} x_i^{e_{i-r}}.
\end{displaymath}

From this we can prove the following theorem, using methods very similar to those found in \cite{FK09}, in which the theorem was proved in the $\Hom{}$-finite case.  Note that the theorem was shown in \cite{DWZ09} using decorated representations and in \cite{Nag09} using Donaldson--Thomas theory.

\begin{theorem}\label{theo::Fconj}
Conjecture 5.6 of \cite{FZ07} holds for $\cA$, that is, any $F$-polynomial has constant term $1$.
\end{theorem}
\demo{ It suffices to show that the polynomial $F_{T_j}$ defined above has constant term $1$.  In order to do so, we will prove that, for any $r< i \leq 2r$, the  number $[\ind{\Gamma} T_j : \Gamma_i]$ vanishes.

We know that $T_j$ lies in the subcategory $\cU$ defined in Section \ref{sect::index g-vectors}, that is, for any $r<i \leq 2r$, the space $\Hom{\cC}(\Sigma^{-1}\Gamma_i, T_j)$ vanishes.  Using \cite[Proposition 2.15]{Plamondon09}, we get that $\Hom{\cC}( T_j, \Sigma\Gamma_i)$ also vanishes.

Let $\overline{T}_1 \rightarrow \overline{T}_0 \rightarrow T_j \rightarrow \Sigma\overline{T}_1$ be a minimal $(\add \Gamma)$-presentation of $T_j$.  Let $r<i\leq 2r$ be a vertex of $Q$.

Suppose that $\Gamma_i$ is a direct summand of $\overline{T}_1$.  Since $\Hom{\cC}( T_j, \Sigma\Gamma_i)$ is zero, this implies that the presentation has the triangle
\begin{displaymath}
	\xymatrix{\Gamma_i\ar[r]^{1_{\Gamma_i}} & \Gamma_i\ar[r] & 0\ar[r] & \Sigma\Gamma_i}
\end{displaymath}
as a direct summand, contradicting the minimality of the presentation.  Thus $\Gamma_i$ is not a direct summand of $\overline{T}_1$.

Suppose that $\Gamma_i$ is a direct summand of $\overline{T}_0$.  Since $i$ is a sink in $Q$, and since $\Gamma_i$ is not a direct summand of $\overline{T}_1$, we get that $\Hom{\cC}(\overline{T}_1, \Gamma_i)$ is zero.  This implies that $\Gamma_i$ is a direct summand of $T_j$, and since the latter is indecomposable, we get that it is isomorphic to the former.  This is a contradiction, since $T$ must be basic.
}

\begin{definition}\label{defi::F}
For any object $M$ of $\cD$, the \emph{$F$-polynomial of $M$} is the polynomial
\begin{displaymath}
	F_M = \sum_{e}\chi \Big( \Gr{e} \big( \Hom{\cC}(\Sigma^{-1}\Gamma, M) \big) \Big) \prod_{i=r+1}^{2r}x_i^{e_{i-r}}
\end{displaymath}
in $\bZ[x_{r+1}, \ldots, x_{2r}]$.
\end{definition}

Thanks to Theorem \ref{theo::Fconj}, this definition is in accordance with the $F_{T_i}$ used above.  Note that we have the equality
\begin{displaymath}
	X'_M \Big|_{x_1 = \ldots = x_r = 1} = \prod_{i=r+1}^{2r} x_i^{[\ind{\Gamma} M : \Gamma_i]} F_M
\end{displaymath}

We can deduce from the multiplication formula of \cite[Proposition 3.16]{Plamondon09} an equality for the polynomials $F_M$.  This was first proved implicitly in \cite[Section 5.1]{Palu08}, see also \cite[Theorem 6.12]{K10}.

\begin{proposition}
Let $M$ and $N$ be objects of $\cD$ such that the space $\Hom{\cC}(M, \Sigma N)$ is one-dimensional.  Let 
\begin{displaymath}
	M \longrightarrow E \longrightarrow N \longrightarrow \Sigma M \quad \textrm{and} \quad N \longrightarrow E' \longrightarrow M \longrightarrow \Sigma N
\end{displaymath}
be non-split triangles.  Then
\begin{displaymath}
	F_M F_N = \prod_{i=r+1}^{2r} x_i^{d_{i-r}} F_E + \prod_{i=r+1}^{2r} x_i^{d'_{i-r}} F_{E'},
\end{displaymath}
where $d = (d_1, \ldots, d_{2r})$ is the dimension vector of the kernel $K$ of the morphism \ $\Hom{\cC}(\Sigma^{-1}\Gamma, M) \longrightarrow \Hom{\cC}(\Sigma^{-1}\Gamma, E)$ and $d' = (d'_1, \ldots, d'_{2r})$ is the dimension vector of the kernel $K'$ of $\Hom{\cC}(\Sigma^{-1}\Gamma, N) \longrightarrow \Hom{\cC}(\Sigma^{-1}\Gamma, E')$. 
\end{proposition}
\demo{ We know from \cite[Proposition 3.16]{Plamondon09} that $X'_M X'_N = X'_E + X'_{E'}$.  Specializing at $x_1 = \ldots = x_r = 1$, we get the equality
\begin{displaymath}
	\prod_{i=r+1}^{2r}x_i^{[\ind{\Gamma}M : \Gamma_i] + [\ind{\Gamma}N : \Gamma_i]} F_M F_N = \prod_{i=r+1}^{2r}x_i^{[\ind{\Gamma}E : \Gamma_i]}F_E + \prod_{i=r+1}^{2r}x_i^{[\ind{\Gamma}E' : \Gamma_i]} F_{E'}.
\end{displaymath}

It follows from \cite[Lemma 3.5]{Plamondon09} (applied to the above triangles shifted by $\Sigma^{-1}$, and with $T= \Sigma^{-1}\Gamma$) that
\begin{eqnarray*}
	\ind{\Gamma}M + \ind{\Gamma}N & = & \ind{\Gamma}E + \ind{\Gamma}K + \ind{\Gamma}\Sigma K \\
	                              & = & \ind{\Gamma}E' + \ind{\Gamma}K' + \ind{\Gamma}\Sigma K',
\end{eqnarray*}
where $K$ and $K'$ are as in the statement of the Proposition.

But $\ind{\Gamma}K + \ind{\Gamma}\Sigma K = \iota(d)$, and using our computation of $\iota(e)$ of Section \ref{sect::index g-vectors}, we get that $-\iota(d)_i = d_{i-r}$ for $r < i \leq 2r$.

Similarly, we get that $\ind{\Gamma}K' + \ind{\Gamma}\Sigma K' = \iota(d')$, and that $-\iota(d')_i = d'_{i-r}$ for $r < i \leq 2r$.

The desired equality follows.

}

%========================================================================================================
\section{Link with decorated representations}\label{sect::dwz}
In this section, an explicit link between cluster categories and the decorated representations of \cite{DWZ08} is established.  We show that the mutation of decorated representations of \cite{DWZ08} corresponds to the derived-equivalence of \cite{KY09}, and we give an interpretation of the $E$-invariant of \cite{DWZ09} as half the dimension of the space of selfextensions of an object in the cluster category.

%--------------------------------------------------------------------------------------
\subsection{Mutations}\label{sect::bijection}
Let $(Q,W)$ be a quiver with potential.  Let $\Gamma = \Gamma_{Q,W}$ be the associated complete Ginzburg dg algebra, and $\cC = \cC_{Q,W}$ be the associated cluster category.  Let $B = B_{Q,W}$ be the endomorphism algebra of $\Gamma$ in $\cC$.  Recall from \cite[Lemma 2.8]{KY09} that $B$ is the Jacobian algebra of $(Q,W)$.  Denote by $F$ the functor $\Hom{\cC}(\Sigma^{-1}\Gamma, ?)$ from $\cC$ to $\Mod B$.  Let $\cD = \cD_{Q,W}$ be the full subcategory of $\pr{\cC}\Gamma \cap \pr{\cC}\Sigma^{-1}\Gamma$ whose objects are those $X$ such that $FX$ is finite dimensional.

Consider the map $\Phi = \Phi_{Q,W}$ from the set of isomorphism classes of objects in $\cD$ to the set of isomorphism classes of decorated representations of $(Q,W)$ defined as follows.  For any object $X$ of $\cD$, write $X = X'\oplus \bigoplus_{i\in Q_0}(e_i \Gamma)^{m_i}$, where $X'$ has no direct summands in $\add \Gamma$.  Such a decomposition of $X$ is unique up to isomorphism, since $\pr{\cC}\Gamma$ is a Krull--Schmidt category, as shown in \cite{Plamondon09}.  Define $\Phi(X)$ to be the decorated representation $(F(X'), \bigoplus_{i\in Q_0}S_i^{m_i} )$, where $(0,S_i)$ is the negative simple representation at the vertex $i$, for any $i$ in $Q_0$.

Consider also the map $\Psi = \Psi_{Q,W}$ from the set of isomorphism classes of decorated representations of $(Q,W)$ to the set of isomorphism classes of objects in $\cD$ defined as follows.  Recall from \cite{Plamondon09} that $F$ induces an equivalence $\pr{\cC}\Sigma^{-1}\Gamma/(\Gamma) \rightarrow \MOD B$, where $\MOD B$ is the category of finitely presented $B$-modules.  Let $G$ be a quasi-inverse equivalence. For any decorated representation $(M, \bigoplus_{i\in Q_0}S_i^{m_i})$, choose a representative $\overline{M}$ of $G(M)$ in $\cD$ which has no direct summands in $\add \Gamma$ (the representative can be chosen to be in $\cD$ thanks to \cite[Lemma 3.2]{Plamondon09}).  Such a representative is unique up to (non-unique) isomorphism.  The map $\Psi$ then sends $(M, \bigoplus_{i\in Q_0}S_i^{m_i})$ to the object $\overline{M} \oplus \bigoplus_{i\in Q_0}(e_i\Gamma_i)^{m_i}$. 

The diagram below summarizes the definitions of $\Phi$ and $\Psi$.
\begin{eqnarray*}
 \left\{ \begin{array}{c}
\textrm{isoclasses of } \\
\textrm{objects of } \cD 
\end{array} \right\}           & \longleftrightarrow & \left\{ \begin{array}{c}
                                                        \textrm{isoclasses of decorated } \\
                                                        \textrm{representations of } (Q,W) 
                                                        \end{array} \right\}  \\
X=X'\oplus\bigoplus_{i=1}^{n}(e_i\Gamma)^{m_i} & \longmapsto  & \Phi(X) = \big(FX', \ \bigoplus_{i=1}^{n}(S_i)^{m_i}\big) \\
\Psi(\mathcal{M}) = \overline{M}\oplus\bigoplus_{i=1}^{n}(e_i\Gamma)^{m_i}& \longmapsfrom & \mathcal{M} = \big( M, \bigoplus_{i=1}^{n}S_i^{m_i}\big)
\end{eqnarray*}

The main result of this subsection states that the maps $\Phi$ and $\Psi$ are mutually inverse bijections, on the one hand, and that, via these maps, the derived equivalences of \cite{KY09} are compatible with the mutations of decorated representations of \cite{DWZ08}, on the other hand.

\begin{proposition}\label{prop::mutations}
With the above notations, $\Phi$ and $\Psi$ are mutual inverse maps.  Moreover, if \ $i\in Q_0$ is not on any cycle of length $\leq 2$, and if $(Q',W') = \widetilde{\mu_i}(Q,W)$, then for any object $X$ of $\cD$, we have that
\begin{displaymath}
	\Phi_{Q',W'}(\widetilde{\mu}_i^-(X)) = \widetilde{\mu}_i(\Phi_{Q,W}(X)),
\end{displaymath}
where the functor $\widetilde{\mu}_i^-$ is as defined after Theorem \ref{theo::KY}
\end{proposition}

The rest of this section is devoted to the proof of the Proposition.

 It is obvious from the definition of $\Phi$ and $\Psi$ that the two maps are mutual inverses.  Thus we only need to show that the two mutations agree.

Let $\Gamma'$ be the complete Ginzburg dg algebra of $(Q',W')$.  Note that $\End{\cC'}(\Gamma')$ is the Jacobian algebra $J(Q',W')$, by \cite[Lemma 2.8]{KY09}.  Let $\cC'$ be the cluster category associated with $(Q',W')$.

We know from \cite{DWZ09} that $\widetilde{\mu}_i(\Phi_{Q,W}(X))$ is a decorated representation of $(Q', W') = \widetilde{\mu}_i(Q,W)$.  We need to show that it is isomorphic to $\Phi_{Q',W'}(\widetilde{\mu}_i^-(X))$.

We can (and will) assume for the rest of the proof that $X$ is indecomposable, as all the maps and functors considered commute with finite direct sums.

We first prove the proposition for some special cases.  

\begin{lemma}
Assume that $X$ is an indecomposable object of $\cD$ such that either
 \begin{itemize}
	\item $X$ is of the form $e_j\Gamma$ for $j\neq i$, or
	\item $X$ is the cone $\Gamma_i^*$ of the morphism 
     \begin{displaymath}
	     \Gamma_i \longrightarrow \bigoplus_{\alpha}\Gamma_{t(\alpha)}
     \end{displaymath}
     whose components are given by left multiplication by $\alpha$.  
 \end{itemize}
Then the equality $\Phi_{Q',W'}(\widetilde{\mu}_i^-(X)) = \widetilde{\mu}_i(\Phi_{Q,W}(X))$ holds.
\end{lemma}
\demo{
Suppose that $X = e_j\Gamma$ for some vertex $i\neq j$.  Then $\widetilde{\mu}_i(\Phi_{Q,W}(X)) = \widetilde{\mu}_i(0, S_j) = (0, S_j)$, and $\Phi_{Q',W'}(\widetilde{\mu}_i^-(X)) = \Phi_{Q',W'}(e_j\Gamma') = (0, S_j)$, so the desired equality holds.

Suppose now that $X$ is the cone $\Gamma_i^*$ of the morphism 
     \begin{displaymath}
	     \Gamma_i \longrightarrow \bigoplus_{\alpha}\Gamma_{t(\alpha)}
     \end{displaymath}
     whose components are given by left multiplication by $\alpha$.  In that case, $\widetilde{\mu}_i^-(X) = e_i\Gamma'$ and $\Phi(X) = (S_i, 0)$, so the desired equality is also satisfied.
}

Now suppose that $X$ is not of the above form.  Using the definition of $\widetilde{\mu}_i^-$, we get that $\Phi(\widetilde{\mu}_i^-(X))$ is equal to $\Phi(\cH om_{\Gamma}(T,X))$, where $T$ is as defined in section \ref{sect::derived}.  Because of our assumptions on $X$, this decorated representation is given by $\Big(\Hom{\cC'}\big(\Sigma^{-1}\Gamma', \cH om_{\Gamma}(T,X)\big), 0\Big)$.

We have the isomorphisms of $\End{\cC'}(\Gamma')$-modules
\begin{eqnarray*}
	\Hom{\cC'}\big(\Sigma^{-1}\Gamma', \cH om_{\Gamma}(T,X)\big) & = & \Hom{\cD\Gamma'}\big(\Sigma^{-1}\Gamma', \cH om_{\Gamma}(T,\overline{X})\big) \\
	 & = & \Hom{\cD\Gamma}(\Sigma^{-1}\Gamma' \otimes^{L}_{\Gamma'} T , \overline{X}) \\
	 & = & \Hom{\cD\Gamma}(\Sigma^{-1}T , \overline{X}) \\
	 & = & \Hom{\cC}(\Sigma^{-1}T , X),
\end{eqnarray*}
where $\overline{X}$ is a lift of $X$ in $\pr{\cD\Gamma}\Sigma^{-1}\Gamma$.

Using this, we prove the Proposition for another special case.

\begin{lemma}
If $X=e_i\Gamma$, then $\Phi_{Q',W'}(\widetilde{\mu}_i^-(X)) = \widetilde{\mu}_i(\Phi_{Q,W}(X))$.
\end{lemma}
\demo{ We have that $\widetilde{\mu}_i(\Phi_{Q,W}(e_i\Gamma)) = (S_i, 0)$.  Moreover, the above calculation gives that $\Phi_{Q',W'}(\widetilde{\mu}_i^-(e_i\Gamma)) = \big(\Hom{\cC}(\Sigma^{-1}T , e_i\Gamma) ,0\big)$.

For any vertex $j\neq i$, we have  $\Hom{\cC}(\Sigma^{-1}T , e_i\Gamma)e_j = \Hom{\cC}(\Sigma^{-1}(e_jT) , e_i\Gamma) = \Hom{\cC}(\Sigma^{-1}(e_j\Gamma) , e_i\Gamma)$, and this space is zero.  

For the vertex $i$, we have that $\Hom{\cC}(\Sigma^{-1}T , e_i\Gamma)e_i = \Hom{\cC}(\Sigma^{-1}(e_iT) , e_i\Gamma) = \Hom{\cC}(\Sigma^{-1}\Gamma_i^* , e_i\Gamma)$, and this space is one-dimensional.

Therefore $\Hom{\cC}(\Sigma^{-1}T , e_i\Gamma)$ is the simple module at the vertex $i$, and this proves the desired equality.
}

We now treat the remaining cases, that is, those where $X$ is not in $\add \Gamma$ and is not $\Gamma_i^*$.  Then $\Phi(X) = (FX, 0)$, and $\widetilde{\mu}_i(\Phi_{Q,W}(X)) = \widetilde{\mu}_i(FX, 0) = (M', 0)$ is computed using section \ref{sect::deco}.  We will show that $\Hom{\cC}(\Sigma^{-1}T , X)$ is isomorphic to $M'$ as a $J(Q',W')$-module, using heavily the definition of $T$ given in section \ref{sect::derived}.

\begin{lemma}\label{lemm::vertices}
For any vertex $j$, the vector spaces $M'e_j$ and $\Hom{\cC}(\Sigma^{-1}T , X)e_j$ are isomorphic.
\end{lemma}
\demo{
If $j$ is a vertex different from $i$, then we have the isomorphisms of vector spaces $\Hom{\cC}(\Sigma^{-1}T , X)e_j = \Hom{\cC}(\Sigma^{-1}(e_jT) , X) = \Hom{\cC}(\Sigma^{-1}(e_j\Gamma) , X) = (FX)e_j = M'e_j$.  

For the vertex $i$, we have that $\Hom{\cC}(\Sigma^{-1}T , X)e_i = \Hom{\cC}(\Sigma^{-1}(e_iT) , X) = \Hom{\cC}(\Sigma^{-1}\Gamma_i^* , X)$.  Let us show that this space is isomorphic to $M'e_i$.

We have triangles in $\cC$
\begin{displaymath}
e_i\Gamma \longrightarrow \bigoplus_{s(a) = i}e_{t(a)}\Gamma \longrightarrow \Gamma_i^* \longrightarrow \Sigma(e_i\Gamma)  \quad \textrm{and}
\end{displaymath}
\begin{displaymath}
 \overline{\Gamma}_i^* \longrightarrow \bigoplus_{t(a) = i}e_{s(a)}\Gamma \longrightarrow e_i\Gamma \longrightarrow \Sigma\Gamma_i^*.
\end{displaymath}

These triangles yield a diagram with exact rows
\begin{displaymath}
	\xymatrix{ (\Sigma^{-1}\overline{\Gamma}_i^*, X)\ar[d]^{\varphi_i^*} & (\Sigma^{-1}\bigoplus_{t(a) = i}e_{s(a)}\Gamma, X)\ar[d]^{-\gamma}\ar[l]_{h\phantom{xxxxx}} & (\Sigma^{-1}(e_i\Gamma), X)\ar@{=}[d]\ar[l]_{\phantom{xxxxx}\beta} & (\Gamma_i^*, X)\ar[l] \\
	(\Sigma^{-1}\Gamma_i^*, X)\ar[r]^{g\phantom{xxxxx}} & (\Sigma^{-1}\bigoplus_{s(a) = i}e_{t(a)}\Gamma, X)\ar[r]^{\phantom{xxxxx}\alpha} & (\Sigma^{-1}(e_i\Gamma), X)\ar[r] & (\Sigma^{-2}\Gamma_i^*, X),
	}
\end{displaymath}
where we write $(Y_1, Y_2)$ for $\Hom{\cC}(Y_1, Y_2)$, where $-\gamma = g\varphi_i^* h$, and where $\varphi_i$ was defined in section \ref{sect::derived}.  Note that $\varphi_i^*$ is an isomorphism.

Notice that, in the notations of section \ref{sect::deco}, $(\Sigma^{-1}\bigoplus_{t(a) = i}e_{s(a)}\Gamma, X) = (FX)_{out}$ and $(\Sigma^{-1}\bigoplus_{s(a) = i}e_{t(a)}\Gamma, X) = (FX)_{in}$.  Moreover, the maps $\alpha$ and $\beta$ in the diagram above correspond to the maps $\alpha$ and $\beta$ of section \ref{sect::deco}.

The map $\gamma$ above also corresponds to the map $\gamma$ defined in section \ref{sect::deco}.  This follows from the computation we made in Remark \ref{rema::W}.

Using the above diagram, we get isomorphisms 
\begin{eqnarray*}
(\Sigma^{-1}\Gamma_i^*, X) & \cong & \Ima g \oplus \Ker g \\
 & \cong & \Ker \alpha \oplus \Ker g
\end{eqnarray*}
and
\begin{eqnarray*}
 \Ker \gamma & \cong & h^{-1}\big(\varphi_i^{*-1}(\Ker g)\big) \\
     & \cong & \Ker h \oplus \Ker g \\
     & \cong & \Ima \beta \oplus \Ker g.
\end{eqnarray*}

Thus $(\Sigma^{-1}\Gamma_i^*, X)$ is (non-canonically) isomorphic to $\Ker\alpha \oplus \frac{\Ker\gamma}{\Ima \beta}$, which is in turn isomorphic to $\frac{\Ker\gamma}{\Ima \beta} \oplus \Ima \gamma \oplus \frac{\Ker\alpha}{\Ima \gamma}$.  But this is precisely $M'e_i$.  
}

It remains to be shown that the action of the arrows of $Q'$ on $\Hom{\cC}(\Sigma^{-1}T , X)$ is the same as on $M'$ in order to get the following Lemma.

\begin{lemma}
As a $J(Q',W')$-module, $\Hom{\cC}(\Sigma^{-1}T, X)$ is isomorphic to $M'$.
\end{lemma}
\demo{  We know from Lemma \ref{lemm::vertices} that the two modules considered are isomorphic as $R$-modules, where $R$ is as in section \ref{sect::QP}.

Now let $a$ be an arrow of $Q$ not incident with $i$.  Then $a$ is an arrow of $Q'$, and its action on $\Hom{\cC}(\Sigma^{-1}T, X)$ is obviously the same as its action on $M'$.

Consider now an arrow of $Q'$ of the form $[ba]$, where $t(a) = i = s(b)$ in $Q$.  By the definition of $M'$ given in section \ref{sect::deco}, $[ba]$ acts as $ba$ on $M'$, that is, $M'_{[ba]} = (FX)_{ba}$.

According to the definition of $T$ given in section \ref{sect::derived}, $[ba]$ acts on $T$ as the map

\begin{eqnarray*}
 T_{s(a)} & \longrightarrow & T_{t(b)} \\
 x & \longmapsto & bax.
\end{eqnarray*}
Hence the action of $[ba]$ on $\Hom{\cC}(\Sigma^{-1}T, X)$ is also given by multiplication by $ba$.  Thus the action of $[ba]$ on $M'$ and on $\Hom{\cC}(\Sigma^{-1}T, X)$ coincide.

There remains to be considered the action of the arrows incident with $i$. 

Keep the notations introduced in the proof of Lemma \ref{lemm::vertices}. We assert that the maps $\varphi_i^* h$ and $g$ encode the action of the arrows incident with $i$.

Recall that in $\cD\Gamma$, the object $\Gamma_i^*$ is isomorphic as a graded module to 
\begin{displaymath}
	\Sigma(e_i\Gamma) \oplus \bigoplus_{\substack{a\in Q_1 \\ s(a) = i}} e_{t(a)}\Gamma,
\end{displaymath}
and that the map $\bigoplus_{\substack{a\in Q_1 \\ s(a) = i}} e_{t(a)}\Gamma \longrightarrow \Gamma_i^*$ is the canonical inclusion.  Thus, its components are given by
\begin{eqnarray*}
 e_{t(a)}\Gamma & \longrightarrow & \Gamma_i^* \\
 x & \longmapsto & e_a x.
\end{eqnarray*}
for any arrow $a$ of $Q$ such that $s(a) = i$.  By the definition of $T$, this is multiplication by $a^{\star}$.  Therefore $g$ encodes the action of the arrows $a^{\star}$ of $Q'$, where $s(a) = i$ in $Q$.

Similarly, recall that in $\cD\Gamma$, the object $\overline{\Gamma}_i^*$ is isomorphic as a graded module to
\begin{displaymath}
	\big( \bigoplus_{\substack{b\in Q_1 \\ t(b) = i}} e_{t(b)}\Gamma \big) \oplus \Sigma^{-1}(e_i\Gamma)
\end{displaymath}
and that the map $\overline{\Gamma}_i^* \longrightarrow \bigoplus_{\substack{b\in Q_1 \\ t(b) = i}} e_{t(b)}\Gamma$ is given by the canonical projection.  Thus its composition with $\varphi_i^*$ is given by the matrix 
$\left( \begin{array}{cc}
-b^{\ast} & -\partial_{ab}W \\
\end{array} \right)$.  Its components are the maps
\begin{eqnarray*}
 \Gamma_i^* & \longrightarrow & e_{s(b)}\Gamma \\
 e_{\Sigma i}x_i + \sum_{s(a) = i}e_ax_a & \longmapsto & -b^{\ast}x_i + \sum_{s(a)=i}(\partial_{ab}W)x_a
\end{eqnarray*}
for any arrow $b$ of $Q$ such that $t(b) = i$.  By the definition of $T$, this is multiplication by $b^{\star}$.  Thus $\varphi^*_i h$ encodes the action of the arrows $b^{\star}$ of $Q'$, where $t(b) = i$ in $Q$.

Finally, recall from Lemma \ref{lemm::vertices} that $\Hom{\cC}(\Sigma^{-1}\Gamma_i^*, X)$ is isomorphic to $\frac{\Ker \gamma}{\Ima \beta}\oplus \Ima \gamma \oplus \frac{\Ker \alpha}{\Ima \gamma}$.  Recall that the summand $\frac{\Ker \gamma}{\Ima \beta}$ corresponds to $\Ker g$, while the summand $\Ima \gamma \oplus \frac{\Ker \alpha}{\Ima \gamma}$ corresponds to $\Ima g$.  

We choose a spliting $\Ima \gamma \oplus \frac{\Ker \alpha}{\Ima \gamma}$ in such a way that $\Ima \varphi_i^* h \cap \frac{\Ker \alpha}{\Ima \gamma} = 0$.  In that case, $g$ is given in matrix form by $\left( \begin{array}{ccc}
0 & \iota & \iota\sigma \\
\end{array} \right)$ and $\varphi_i^* h$, by $\left( \begin{array}{c}
-\pi\rho \\
-\gamma \\
0
\end{array} \right)$, in the notations of section \ref{sect::deco}.

This proves that the action of the arrows of $Q'$ on $M'$ and on $\Hom{\cC}(\Sigma^{-1}T, X)$ coincide, finishing the proof of the Lemma.
}

We have proved Proposition \ref{prop::mutations}.

%--------------------------------------------------------------------------------------
\subsection{Interpretation of $F$-polynomials, $\bg$-vectors and $h$-vectors}\label{sect::fg}
In this section, we study the relation between the $F$-polynomials of objects of $\cD$ and of decorated representations, and between the index of objects in $\cD$ and the $\bg$-vectors of decorated representations.  We also give an interpretation of the $h$-vector.

Let $(Q,W)$ be a quiver with potential, and let $\cC$ be the associated cluster category.  We keep the notations of the previous section for the maps $\Phi$ and $\Psi$.

We first prove a result regarding $F$-polynomials.

\begin{proposition}\label{prop::2 Fpolynomials}
Let $X$ be an object of $\cD$.  Then we have the equality
 \begin{displaymath}
   F_X(x_{r+1}, \ldots, x_n) = F_{\Phi(X)}(x_{r+1}, \ldots, x_{n}).
 \end{displaymath}
\end{proposition}
\demo{  This is immediate from the definitions of $F_X$, $\Phi$ and $F_{\Phi(X)}$, given in Definition \ref{defi::F}, Section \ref{sect::bijection} and Section \ref{sect::invariants}, respectively.
}

We now prove that $\bg$-vectors of decorated representations and indices of objects in the cluster category are closely related.  

We will need the following $\Hom{}$-infinite extension of \cite[Lemma 7]{Palu08}.
\begin{lemma}\label{lemm::indices ext}
Let $M$ be an indecomposable object of $\cD$.  Then
\begin{displaymath}
	[\ind{\Gamma}M : e_i\Gamma] = \left\{ \begin{array}{ll}
\delta_{ij} & \textrm{if $M \cong e_i\Gamma$}\\
\dim \Ext{1}{B}(S_i, FM) - \dim \Hom{B}(S_i, FM) & \textrm{otherwise},
\end{array} \right.
\end{displaymath}
where $B=\End{\cC}(\Gamma)$.
\end{lemma}
\demo{The result is obvious if $M$ lies in $\add\Gamma$.  Suppose it does not.  Let $T_1 \longrightarrow T_0 \longrightarrow M \longrightarrow \Sigma T_1$ be an $(\add\Gamma)$-presentation of $M$.  

The opposite category $\cC^{op}$ is triangulated, with suspension functor $\Sigma_{op} = \Sigma^{-1}$.  Thus, in $\cC^{op}$, we have a triangle $\Sigma^{-1}_{op} T_0 \longrightarrow \Sigma^{-1}_{op} T_1 \longrightarrow M \longrightarrow  T_0$.  Applying the functor $F' = \Hom{\cC^{op}}(\Sigma^{-1}_{op}\Gamma, ?)$, we get a minimal projective resolution
\begin{displaymath}
	F'\Sigma^{-1}_{op}T_0 \longrightarrow F'\Sigma^{-1}_{op} T_1 \longrightarrow F'M \longrightarrow 0
\end{displaymath}
of $F'M$ as a $B^{op}$-module.

Letting $S'_i$ be the simple $B^{op}$-module at the vertex $i$, we apply $\Hom{B'}(?, S_i)$ to the above exact sequence and get a complex
\begin{displaymath}
	0 \longrightarrow \Hom{B^{op}}(F'\Sigma_{op}^{-1}T_1, S'_i) \longrightarrow \Hom{B^{op}}(F'\Sigma^{-1}_{op}T_0, S'_i) \longrightarrow \ldots
\end{displaymath}
whose differential vanishes, since the presentation is minimal.

Therefore we have the equalities
\begin{eqnarray*}
  [\ind{\Gamma}M : e_i\Gamma]  & = & \dim \Ext{1}{B^{op}}(F'M, S_i') - \dim \Hom{B^{op}}(F'M, S_i') \\
    & = & \dim \Ext{1}{B}(S_i, DF'M) - \dim \Hom{B}(S_i, DF'M),
\end{eqnarray*}
where $S_i$ is the simple $B$-module at the vertex $i$.

Now, using \cite[Proposition 2.16]{Plamondon09}, we get that
\begin{displaymath}
	DF'M = D\Hom{\cC^{op}}(\Sigma^{-1}_{op}\Gamma, M) = D\Hom{\cC}(M, \Sigma\Gamma) \cong \Hom{\cC}(\Sigma^{-1}\Gamma, M) = FM.
\end{displaymath}
Thus $DF'M$ is isomorphic to $FM$ as a $B$-module.  This proves the lemma.
}

We now prove the result on $\bg$-vectors of decorated representations.

\begin{proposition}\label{prop::2 gvectors}
Let $(Q,W)$ be a quiver with potential, and let $\cC$ be the associated cluster category.  Let $X$ be an object of $\cD$. Let $\bg_{\Phi(X)} = (g_1, \ldots, g_n)$ be the $\bg$-vector of the decorated representation $\Phi(X)$.  Then we have the equality
  \begin{displaymath}
	  g_i = [\ind{\Gamma}X : \Gamma_i]
  \end{displaymath}
for any vertex $i$ of $Q$.
\end{proposition}
\demo{ We can assume that $X$ is indecomposable.  If $X$ lies in $\add\Gamma$, then the result is obviously true.  Suppose that $X$ does not lie in $\add\Gamma$.

Using the two triangles in $\cC$
\begin{displaymath}
e_i\Gamma \longrightarrow \bigoplus_{s(a) = i}e_{t(a)}\Gamma \longrightarrow \Gamma_i^* \longrightarrow \Sigma(e_i\Gamma)  \quad \textrm{and}
\end{displaymath}
\begin{displaymath}
 \overline{\Gamma}_i^* \longrightarrow \bigoplus_{t(a) = i}e_{s(a)}\Gamma \longrightarrow e_i\Gamma \longrightarrow \Sigma\Gamma_i^*.
\end{displaymath}
and applying the functor $F = \Hom{\cC}(\Sigma^{-1}\Gamma, ?)$, we get a projective resolution of the simple $B$-module $S_i$ at the vertex $i$:
\begin{displaymath}
	P_i \longrightarrow \bigoplus_{s(a) = i}P_{t(a)} \longrightarrow \bigoplus_{t(a) = i}P_{s(a)} \longrightarrow P_i \longrightarrow S_i \longrightarrow 0,
\end{displaymath}
where $P_j$ is the indecomposable projective $B$-module at the vertex $j$.  Applying now the functor $\Hom{B}(?, FM)$, we get the complex
\begin{displaymath}
	0 \longrightarrow (FM)_i \stackrel{\beta_i}{\longrightarrow} (FM)_{out} \stackrel{-\gamma_i}{\longrightarrow} (FM)_{in} \stackrel{\alpha_i}{\longrightarrow} (FM)_i.
\end{displaymath}

From this complex, we see that $\Hom{B}(S_i, M) = \Ker \beta_i$ and that $\Ext{1}{B}(S_i, M) = \Ker\gamma_i / \Ima \beta_i$.  We also deduce an exact sequence
\begin{displaymath}
	0 \longrightarrow \Ker \beta_i \longrightarrow (FM)_i \stackrel{\beta_i}{\longrightarrow} \Ker \gamma_i \longrightarrow \Ker\gamma_i / \Ima \beta_i \longrightarrow 0.
\end{displaymath}
Using the above arguments and Lemma \ref{lemm::indices ext}, we get the equalities
\begin{eqnarray*}
	 [\ind{\Gamma}X:e_i\Gamma] & = & \dim \Ext{1}{B}(S_i, M) - \dim \Hom{B}(S_i, M) \\
	  & = &  \dim (\Ker\gamma_i / \Ima \beta_i) \dim \Ker \beta_i \\
	  & = & \dim \Ker \gamma_i - \dim (FM)_i \\
	  & = & g_i. 
\end{eqnarray*}
This finishes the proof.
}

As a corollary of the proof of the above Proposition, we get an interpretation of the $h$-vector of a decorated representation.
\begin{corollary}\label{coro::hvector}
For any decorated representation $\cM = (M,V)$ of a quiver with potential $(Q,W)$, we have the equality
\begin{displaymath}
	h_i = -\dim \Hom{J(Q,W)}(S_i, M)
\end{displaymath}
for any vertex $i$ of $Q$.
\end{corollary}

This provides us with a way of ``counting'' the number of terms in a minimal presentation.
\begin{corollary}\label{coro::counting}
If $\bg = (g_1, \ldots, g_n)$ and $\bh = (h_1, \ldots, h_n)$ are the $\bg$-vector and $\bh$-vector of a decorated representation $\cM = (M,V)$, $\bh' = (h'_1, \ldots, h'_n)$ is the $\bh$-vector of $\mu_i(\cM)$, and if
\begin{displaymath}
	T_1 \longrightarrow T_0 \longrightarrow \Psi(\cM) \longrightarrow \Sigma T_1
\end{displaymath}
is a minimal $(\add \Gamma)$-presentation of $\Psi(\cM)$ (see Proposition \ref{prop::mutations}), then $-h_i$ and $-h'_i$ are the number of direct summands of $T_1$ and $T_0$ which are isomorphic to $\Gamma_i$, respectively.
\end{corollary}
\demo{ It follows from Corollary \ref{coro::hvector} that $-h_i = \dim \Hom{J(Q,W)}(S_i, M)$.

Let $T_i^*$ be an indecomposable object of $\cD$ such that $\Hom{\cC}(\Sigma^{-1}\Gamma, T_i^*)$ is the simple $S_i$.  Then, by \cite[Lemma 3.2]{Plamondon09}, we have that
\begin{displaymath}
	\Hom{J(Q,W)}(S_i, M) \cong \Hom{\cC}(T_i^*, \Psi(M))/(\Gamma).
\end{displaymath}
Applying $\Hom{\cC}(T_i^*, ?)$ to the presentation, we get a long exact sequence
\begin{displaymath}
	(T_i^*, T_0) \stackrel{\psi_*}{\longrightarrow} (T_i^*, \Psi(\cM)) \stackrel{\phi_*}{\longrightarrow} (T_i^*, \Sigma T_1) \longrightarrow (T_i^*, \Sigma T_0).
\end{displaymath}

We see that the image of $\psi_*$ is $(\Gamma)(T_i^*, \Psi(M))$, so that $\Hom{\cC}(T_i^*, \Psi(M))/(\Gamma)$ is isomorphic to the image of $\phi_*$.  Thus $-h_i$ is the dimension of the image of $\phi_*$.

Using \cite[Proposition 2.16]{Plamondon09}, we get that the morphism $(T_i^*, \Sigma T_1) \longrightarrow (T_i^*, \Sigma T_0)$ is isomorphic to the morphism $D(\Sigma^{-1} T_1, T_i^*) \longrightarrow D(\Sigma^{-1} T_0, T_i^*)$, and this morphism is zero since the presentation is minimal.  Thus $\phi_*$ is surjective.

Therefore $-h_i$ is equal to the dimension of $\Hom{\cC}(\Sigma^{-1}T_1, T_i^*)$, which is equal to the number of direct factors of $T_1$ isomorphic to $\Gamma_i$ in any decomposition of $T_1$.

Furthermore, \cite[Lemma 5.2]{DWZ09} gives that $g_i = h_i - h'_i$, and by Proposition \ref{prop::2 gvectors}, $g_i = [\ind{\Gamma}\Psi(\cM):\Gamma_i]$.  This immediately implies that $-h'_i$ is equal to the number of direct factors of $T_0$ isomorphic to $\Gamma_i$, and finishes the proof. 
}

\begin{remark}\label{rema::reformulation}
Corollary \ref{coro::counting} allows us to reformulate Remark \ref{rema::indices} in the following way.  If $M$ is any object of $\cD$, and if $\bh = (h_1, \ldots, h_n)$ and $\bh' = (h'_1, \ldots, h'_n)$ are the $\bh$-vectors of $\Phi(M)$ and $\tilde{\mu}_i\Phi(M)$, respectively, then
\begin{displaymath}
	  [\ind{T'}M : T'_j] = \left\{ \begin{array}{ll}
                       -[\ind{T}M : T_{i}] \qquad (\textrm{if } i=j) \\
                        \phantom{}[\ind{T}M : T_{j}] -h'_i[b_{ji}]_+ + h_i[-b_{ji}]_+ \qquad (\textrm{if } i\neq j).
                   \end{array} \right.
\end{displaymath}
\end{remark}

As a corollary, we get a proof of Conjecture 6.10 of \cite{FZ07}.

\begin{corollary}\label{coro::6.10}
Conjecture 6.10 of \cite{FZ07} is true, that is, if $\bg = (g_1, \ldots, g_n)$ and $\bg' = (g'_1, \ldots, g'_n)$ are the $\bg$-vectors of one cluster variable with respect to two clusters $t$ and $t'$ related by one mutation at vertex $i$, and if $\bh = (h_1, \ldots, h_n)$ and $\bh' = (h'_1, \ldots, h'_n)$ are its $\bh$-vectors with respect to those clusters, then we have that
\begin{displaymath}
	h'_i = -[g_i]_+ \quad \textrm{and} \quad h_i = \min(0,g_i).
\end{displaymath}
\end{corollary}
\demo{ Let $M$ be an indecomposable object of $\cD$ such that $X'_M$ is the cluster variable considered in the statement.  In particular, $M$ is reachable, and thus rigid.  It follows from equation (5.5) of \cite{DWZ09} that the $\bh$-vector of the cluster variable corresponds to the $\bh$-vector of the associated decorated representation.

Since $M$ is rigid, Proposition \ref{lemm::repetition} tells us that any minimal $(\add\Gamma)$-presentation of $M$ has disjoint direct factors.  The result follows directly from this and from Corollary \ref{coro::counting}.
}

\begin{remark}
Conjecture 6.10 of \cite{FZ07} also follows directly from Conjecture 7.12 (see Theorem \ref{theo::gvectors}(4) above) and equations (6.15) and (6.26) of \cite{FZ07}.  We give the above proof because it is an application of the results developped in this paper.
\end{remark}

Finally, we get an interpretation of the substitution formula of \cite[Lemma 5.2]{DWZ09} in terms of cluster characters.

\begin{corollary}\label{coro::subst}
Let $(Q,W)$ be a quiver with potential.  Let $i$ be an admissible vertex of $Q$, and let $\varphi_X : \bQ(x_1', \ldots, x_n') \longrightarrow \bQ(x_1, \ldots, x_n)$ be the field isomorphism sending $x_j'$ to $x_j$ if $i\neq j$ and to
\begin{displaymath}
	(x_i)^{-1}(\prod_{\ell=1}^n x_\ell^{[b_{\ell i}]_+} + \prod_{\ell=1}^n x_\ell^{[-b_{\ell i}]_+})
\end{displaymath}
if $i=j$.  Let $\cC$ and $\cC'$ be the cluster categories of $(Q,W)$ and $\widetilde{\mu}_i(Q,W)$, respectively, and let $\widetilde{\mu}_i^+:\cC' \longrightarrow \cC$ be the associated functor (see \cite[Theorem 3.2]{KY09}).

Then for any object $M$ of the subcategory $\cD'$ of $\cC'$, we have that 
\begin{displaymath}
X'_{\widetilde{\mu}_i^+(M)} = \varphi_X(X'_M).
\end{displaymath}
\end{corollary}
\demo{ Consider the field isomorphism $\varphi_Y:\bQ(y_1', \ldots, y_n') \longrightarrow \bQ(y_1, \ldots, y_n)$ whose action on $y_j'$ is given by
\begin{displaymath}
	       \varphi_Y(y'_j) = \left\{ \begin{array}{ll}
              y_i^{-1} & \textrm{if $i=j$}\\
              y_j y_i^{m}(y_i + 1)^{-m} & \textrm{if there are $m$ arrows from $i$ to $j$}\\
              y_j (y_i + 1)^{m} & \textrm{if there are $m$ arrows from $j$ to $i$.}
                     \end{array} \right.
     \end{displaymath}
Consider also the morphism $\hat{(-)}:\bQ(y_1, \ldots, y_n) \longrightarrow \bQ(x_1, \ldots, x_n)$ sending each $y_j$ to 
\begin{displaymath}
	\hat{y}_j = \prod_{\ell = 1}^n x_\ell^{b_{\ell j}}.
\end{displaymath}
Denote by the same symbol the corresponding map from the field $\bQ(y_1', \ldots, y_n')$ to $\bQ(x_1', \ldots, x_n')$.  Then \cite[Proposition 3.9]{FZ07} implies that $\varphi_X(\hat{z}) = \widehat{(\varphi_Y(z))}$ for any $z\in \bQ(y_1', \ldots, y_n')$.  In other words, the following diagram commutes:
\begin{displaymath}
	\xymatrix{ \bQ(y_1, \ldots, y_n)\ar[r]^{\hat{(-)}} & \bQ(x_1, \ldots, x_n) \\
	           \bQ(y'_1, \ldots, y'_n)\ar[r]^{\hat{(-)}}\ar[u]^{\varphi_Y} & \bQ(x'_1, \ldots, x'_n)\ar[u]^{\varphi_X}.
	}
\end{displaymath}

Let us now compute $\varphi_X(X'_M)$.  We have that
\begin{eqnarray*}
\varphi_X(X'_M) & = & \varphi_X(x'^{\ind{\Gamma'}M}F_M(\hat{y}'_1, \ldots, \hat{y}'_n) \\
 & = & \varphi_X(x'^{\ind{\Gamma'}M}) F_M(\widehat{(\varphi_Y(y'_1))}, \ldots, \widehat{(\varphi_Y(y'_n))}).
\end{eqnarray*}

Now, using \cite[Lemma 5.2]{DWZ09}, we can express the right-hand side of the equation in terms of the $\hat{y}_j$.  The equalities thus continue:
\begin{eqnarray*}
\varphi_X(X'_M) & = & \varphi_X(x'^{\ind{\Gamma'}M}) \varphi_X(1 + \hat{y}'_i)^{-h'_i} (1+\hat{y}_i)^{h_i} F_M(\hat{y}_1, \ldots, \hat{y}_n) \\
 & = & \varphi_X(x'^{\ind{\Gamma'}M}) \varphi_X(1 + \hat{y}'_i)^{-h'_i} (1+\hat{y}_i)^{h_i} x^{-\ind{\Gamma}\widetilde{\mu}_i^+(M)} X'_{\widetilde{\mu}_i^+(M)}.
\end{eqnarray*}

Thus, in order to prove the Corollary, we must show that 
\begin{eqnarray}
	\varphi_X(x'^{\ind{\Gamma'}M}) \varphi_X(1 + \hat{y}'_i)^{-h'_i} (1+\hat{y}_i)^{h_i} x^{-\ind{\Gamma}\widetilde{\mu}_i^+(M)} & = & 1.
	\label{eqn}
\end{eqnarray}

We do this in several steps.  First, using the definition of $\varphi_X$ and $\varphi_Y$, we get
\begin{eqnarray*}
  \varphi_X(1 + \hat{y}'_i)^{-h'_i} (1+\hat{y}_i)^{h_i} & = & (1 + \widehat{\varphi_Y(y'_i)})^{-h'_i} (1+\hat{y}_i)^{h_i} \\
  & = & (1 + \hat{y}_i^{-1})^{-h'_i} (1+\hat{y}_i)^{h_i} \\
  & = & \hat{y}_i^{h'_i}(1+\hat{y}_i)^{h_i - h'_i}.
\end{eqnarray*}

Now, using Proposition \ref{prop::2 gvectors}, we get the equalities
\begin{eqnarray*}
\varphi_X(x'^{\ind{\Gamma'}M})x^{-\ind{\Gamma}\widetilde{\mu}_i^+(M)} & = & \varphi_X(\prod_{\ell = 1}^{n} (x'_\ell)^{g'_\ell}) \prod_{\ell=1}^{n}x_\ell^{-g_\ell}\\
  & = & x_i^{g_i}  (\prod_{\ell=1}^n x_\ell^{[b_{\ell i}]_+} + \prod_{\ell=1}^n x_\ell^{[-b_{\ell i}]_+})^{-g_i} (\prod_{\ell\neq i}x_\ell^{g'_\ell-g_\ell})   x_i^{-g_i} \\
  & = & (\prod_{\ell=1}^n x_\ell^{[b_{\ell i}]_+} + \prod_{\ell=1}^n x_\ell^{[-b_{\ell i}]_+})^{-g_i} (\prod_{\ell\neq i}x_\ell^{g'_\ell-g_\ell}).    
\end{eqnarray*}

Thus we have, using the fact that $g_i = h_i - h'_i$ \cite[Lemma 5.2]{DWZ09}, that the left-hand side of equation (\ref{eqn}) is equal to
\begin{displaymath}
\hat{y}_i^{h'_i}(1+\hat{y}_i)^{g_i}(\prod_{\ell=1}^n x_\ell^{[b_{\ell i}]_+} + \prod_{\ell=1}^n x_\ell^{[-b_{\ell i}]_+})^{-g_i} (\prod_{\ell\neq i}x_\ell^{g'_\ell-g_\ell})      
\end{displaymath}
which is in turn equal to (using Remark \ref{rema::reformulation})
\begin{eqnarray*}
 \hat{y}_i^{h'_i}(\prod_{\ell = 1}^n x_\ell^{-[-b_{\ell i}]_+})^{g_i}(\prod_{\ell\neq i}x_\ell^{g'_\ell-g_\ell}) 
   & = & \hat{y}_i^{h'_i}(\prod_{\ell = 1}^n x_\ell^{-[-b_{\ell i}]_+})^{g_i} (\prod_{\ell\neq i}x_\ell^{h_i[-b_{\ell i}]_+ -h'_i[b_{\ell i}]_+}) \\
   & = & \hat{y}_i^{h'_i}(\prod_{\ell\neq i}x_\ell^{h'_i[-b_{\ell i}]_+ - h_i[-b_{\ell i}]_+  +    h_i[-b_{\ell i}]_+ -h'_i[b_{\ell i}]_+}) \\
   & = & \hat{y}_i^{h'_i}(\prod_{\ell\neq i}x_\ell^{ -h'_i b_{\ell i}}) \\
   & = & \prod_{\ell\neq i}x_\ell^{ h'_i b_{\ell i}  -h'_i b_{\ell i}} \\
   & = & 1.
\end{eqnarray*}

This finishes the proof.
}

%--------------------------------------------------------------------------------------
\subsection{Extensions and the $E$-invariant}\label{sect::E}

In this section, we give an interpretation of the $E$-invariant of a decorated representation, as defined in \cite{DWZ09} (its definition was recalled in section \ref{sect::invariants}), as the dimension of a space of extensions, using the map $\Phi$ of section \ref{sect::bijection}.

\begin{proposition}\label{prop::E}
Let $(Q,W)$ be a quiver with potential, and let $\cC$ be the associated cluster category.  Let $X$ and $Y$ be objects of $\cD$.  Then we have the following equalities:
 \begin{enumerate}
   \item $E^{inj}(\Phi(X), \Phi(Y)) = \dim (\Sigma\Gamma)(X,  \Sigma Y)$;
   \item $E^{sym}(\Phi(X), \Phi(Y)) = \dim (\Sigma\Gamma)(X,  \Sigma Y) +  \dim (\Sigma\Gamma)(Y,  \Sigma X)$;
   \item $E(\Phi(X)) = (1/2)\dim \Hom{\cC}(X,  \Sigma X)$,
 \end{enumerate}
where $(\Sigma\Gamma)(X, Y)$ is the subspace of $\Hom{\cC}(X,Y)$ containing all morphisms factoring through an object of $\add \Sigma\Gamma$.
\end{proposition}
\demo{The second equality follows immediately from the first one.  

The third equality follows from the second one.  Indeed, the second equality implies that $(\Sigma\Gamma)(X, \Sigma X)$ is finite-dimensional.  It then follows from \cite[Lemma 3.8]{Plamondon09} that we have an isomorphism
  \begin{displaymath}
	  (\Sigma\Gamma)(X, \Sigma X) \cong D\Hom{\cC}(X, \Sigma X)/(\Sigma\Gamma).
  \end{displaymath}
Since $\dim \Hom{\cC}(X, \Sigma X) = \dim (\Sigma\Gamma)(X, \Sigma X) + \dim \Hom{\cC}(X, \Sigma X)/(\Sigma\Gamma)$, we get that 
  \begin{eqnarray*}
	  \dim \Hom{\cC}(X, \Sigma X) & = & 2\dim (\Sigma\Gamma)(X, \Sigma X) \\
	         & = & E^{sym}(\Phi(X), \Phi(X)) \\
	         & = & 2 E(\Phi(X)).
  \end{eqnarray*}
  
Let us now prove the first equality.  Let
\begin{displaymath}
	T_1^Y \longrightarrow T_0^Y \longrightarrow Y \longrightarrow \Sigma T_1^Y
\end{displaymath}
be an $(\add \Gamma)$-presentation of $Y$.  This triangle yields an exact sequence
\begin{displaymath}
	\xymatrix{ (X, Y) \ar[r]^{u\phantom{xx}} & (X, \Sigma T_1^Y)\ar[r] & (X, \Sigma T_0^Y)\ar[r] & (X, \Sigma Y)\ar[r]^{v\phantom{xx}} & (X, \Sigma^2 T_1^Y),
	}
\end{displaymath}
which in turn gives an exact sequence
\begin{displaymath}
	\xymatrix{ 0\ar[r] & \Ima u \ar[r] & (X, \Sigma T_1^Y)\ar[r] & (X, \Sigma T_0^Y)\ar[r] &  \Ker v\ar[r] & 0.
	}
\end{displaymath}
Since $X$ is in $\cD$, the two middle terms of this exact sequence are isomorphic to $(T_i^Y, \Sigma X)$ (for $i=1,2$) thanks to \cite[Proposition 2.16]{Plamondon09}, and these spaces are finite-dimensional.  Therefore all of the terms of the exact sequence are finite-dimensional. 

Now, $\Ima u$ is isomorphic to $(X, Y)/\Ker u$, and $\Ker u$ is exactly $(\Gamma)(X, Y)$.  Therefore, by \cite[Lemma 3.2]{Plamondon09}, $\Ima u$ is isomorphic to the space $\Hom{J(Q,W)}(FX, FY)$, where $F = \Hom{\cC}(\Sigma^{-1}\Gamma, ?)$.

Moreover, $\Ker v$ is exactly $(\Sigma\Gamma)(X, \Sigma Y)$.

Thus, using the above exact sequence and Proposition \ref{prop::2 gvectors}, we have the equalities
\begin{eqnarray*}
	\dim (\Sigma\Gamma)(X, \Sigma Y) & = & \dim \Hom{J(Q,W)}(FX, FY) - \dim (X, \Sigma T_1^Y) + \dim (X, \Sigma T_0^Y) \\
	 & = & \dim \Hom{J(Q,W)}(FX, FY) - \dim (T_1^Y, \Sigma X) + \dim (T_0^Y, \Sigma X) \\
	 & = & \dim \Hom{J(Q,W)}(FX, FY) - \sum_{i=1}^n [T_1^Y:T_i](\dim (FX)_i) + \\
	   & & \ + \sum_{i=1}^n [T_0^X:T_i](\dim (FX)_i) \\
	 & = & \dim \Hom{J(Q,W)}(FX, FY) + \sum_{i=1}^n [\ind{\Gamma}\Sigma Y : \Gamma_i](\dim (FX)_i) \\
	 & = & \dim \Hom{J(Q,W)}(FX, FY) + \sum_{i=1}^n \bg_i(\Phi(Y))(\dim (FX)_i) \\
	 & = & E^{inj}(\Phi(X), \Phi(Y)),  
\end{eqnarray*}
where $[T_j^Y : T_i]$ is the number of direct summands of $T_j^Y$ isomorphic to $T_i$ in any decomposition of $T_j^Y$ into indecomposable objects, and where the $\bg$-vector of $\Phi(Y)$ is given by $(\bg_1(\Phi(Y)), \ldots, \bg_n(\Phi(Y)))$.  This finishes the proof.
}

As a corollary, we get the following stronger version of \cite[Lemma 9.2]{DWZ09}.

\begin{corollary}\label{coro::dwz}
Let $\cM$ and $\cM'$ be two decorated representations of a quiver with potential $(Q,W)$.  Assume that $E(\cM') = 0$.  Then the following conditions are equivalent:
  \begin{enumerate}
    \item $\cM$ and $\cM'$ are isomorphic;
    \item $E(\cM) = 0$, and $\bg_{\cM} = \bg_{\cM'}$.
  \end{enumerate}
\end{corollary}
\demo{  Condition (1) obviously implies condition (2).  Now assume that condition (2) is satisfied.  Then Proposition \ref{prop::E} implies that $\Psi(\cM)$ and $\Psi(\cM')$ are rigid objects of $\cD$.  By Proposition \ref{prop::2 gvectors}, the indices of $\Psi(\cM)$ and $\Psi(\cM')$ are given by $\bg_{\cM}$ and $\bg_{\cM'}$.  By hypothesis, their indices are the same.  Thus, by Proposition \ref{prop::rigidindex}, $\Psi(\cM)$ and $\Psi(\cM')$ are isomorphic, and so are $\cM$ and $\cM'$.
}

%-------------------------------------------------------------------------------
{\footnotesize
{}  
}

\end{document}